\documentclass[preprint,3p,12pt]{elsarticle}
\usepackage{amsmath}
\usepackage{stmaryrd}
\usepackage{bbding}
\usepackage{dcolumn}
\usepackage{graphicx}
\usepackage{amsfonts}
\usepackage{amssymb}
\usepackage{psfrag}
\usepackage{wrapfig}
\usepackage{subfigure}
\usepackage{makeidx}
\usepackage{bm}
\usepackage{epsf}
\usepackage{epsfig}
\usepackage{setspace}
\usepackage{color}
\usepackage{booktabs}
\usepackage{setspace}
\usepackage[ruled,linesnumbered]{algorithm2e}

\definecolor{SpringGreen4}{RGB}{0,139,69}

\begin{document}
% \title{Compact gas-kinetic scheme based high-resolution scale-resolving simulation for compressible flows on structured meshes}
% \title{Efficient implementation of high-order compact gas-kinetic scheme for compressible flows on structured meshes}
\title{An effective implementation of high-order compact gas-kinetic scheme on structured meshes for compressible flows}

\author[HKUST1]{Yaqing Yang}
\ead{yangyq@ust.hk}  %yqyangust@163.com

\author[HKUST1]{Fengxiang Zhao}
\ead{fzhaoac@connect.ust.hk}

\author[HKUST1,HKUST2]{Kun Xu\corref{cor}}
\ead{makxu@ust.hk}

\cortext[cor]{Corresponding author}
\address[HKUST1]{Department of Mathematics, Hong Kong University of Science and Technology, Clear Water Bay, Kowloon, Hong Kong}
\address[HKUST2]{Shenzhen Research Institute, Hong Kong University of Science and Technology, Shenzhen, China}
 
\begin{abstract}

A novel fifth-order compact gas-kinetic scheme is developed for high-resolution simulation of compressible flows on structured meshes. 
Its accuracy relies on a new multidimensional fifth-order compact reconstruction that uses line-averaged derivatives to introduce additional degrees of freedom, enabling a compact stencil with superior resolution. 
For non-orthogonal meshes, reconstruction is performed on a standard reference cell in a transformed computational space. 
This approach provides a unified polynomial form, significantly reducing memory usage and computational cost while simplifying implementation compared to direct multi-dimensional or dimension-by-dimension methods. 
A nonlinear adaptive method ensures high accuracy and robustness by smoothly transitioning from the high-order linear scheme in smooth regions to a second-order scheme at discontinuities. 
The method is implemented with multi-GPU parallelization using CUDA and MPI for large-scale applications. 
Comprehensive numerical tests, from subsonic to supersonic turbulence, validate the scheme's high accuracy, resolution and excellent robustness.

\end{abstract}
 

\maketitle

\section{Introduction}

Scale-resolving simulations, including large eddy simulation (LES) and direct numerical simulation (DNS), aim to explicitly resolve a range of turbulent scales. DNS, in particular, resolves the entire spectrum of turbulent scales, making it an indispensable tool for understanding turbulence mechanisms.
The exceptional solution accuracy and computational efficiency of high-order methods have led to their widespread adoption in LES and DNS.
Extensive efforts have been dedicated to developing high-order schemes, including the high-order finite difference (FD) methods \cite{FD-1,FD-2}, the essential non-oscillatory (ENO) scheme \cite{ENO-1,ENO-2}, weighted essential non-oscillatory (WENO) scheme \cite{WENO-1,WENO-2}, Hermite WENO (HWENO) scheme \cite{HWENO-1}, discontinuous Galerkin (DG) method \cite{DG-0,DG-1}, flux reconstruction (FR) \cite{FR}, correction procedure using reconstruction (CPR) method \cite{CPR}, among others.
Finite difference and finite volume methods typically achieve high-order accuracy by using expanded stencils that incorporate additional neighboring information for reconstruction.
In contrast, DG-type methods improve accuracy by increasing the degree of interpolation polynomials within each cell, thereby enhancing the internal degrees of freedom and naturally exhibiting compactness. 
However, DG-type methods inevitably require identifying troubled cells and applying additional limiting procedures to handle cases with discontinuities \cite{DG-1,DDG}.

In recent years, high-order compact gas-kinetic schemes (GKS) have been developed for the high-fidelity simulation of compressible flows \cite{CGKS-high-2,CGKS-high-3,CGKS-high-4}. These schemes are founded upon the gas distribution function from the GKS framework \cite{GKS-Xu1,GKS-Xu2,GKS-Xu3}, which offers several intrinsic advantages: the simultaneous computation of both inviscid and viscous fluxes from a single unified formulation; an adaptive mechanism that transitions between equilibrium state fluxes in smooth regions and non-equilibrium transport fluxes suitable for discontinuities; and a single-step, second-order time-marching capability.
Building upon this foundation, the high-order compact GKS further leverages the gas distribution function to model the time-accurate evolution of flow variables at cell interfaces. This process, in turn, provides a direct update for the cell-averaged derivatives of conservative variables, which forms the basis for the compact high-order spatial reconstruction. Furthermore, capitalizing on the time-accurate nature of the gas distribution function, an efficient two-stage, fourth-order (S2O4) temporal discretization method \cite{s2o4-0,GRP-high-1,s2o4-1} has also been utilized.
Consequently, high-order compact GKS maintains robustness comparable to second-order schemes for problems with strong shocks \cite{CGKS-high-6}, as they update the derivatives within the cell through the evolved solution at the cell interface, eliminating the need for dedicated shock detection and a posteriori modification procedures.
Previous studies on CGKS have primarily focused on two-dimensional and three-dimensional unstructured meshes \cite{CGKS-high-3,CGKS-high-6,CGKS-high-4}. 
For structured meshes, a higher-order and high-resolution GKS was proposed in \cite{CGKS-high-2}, demonstrating spectral-like resolution characteristics in one dimension. 
This highlights the potential of compact GKS for resolving the multi-scale flow structures in turbulence with high fidelity. 
However, its extension to high-dimensional cases relies on a dimension-by-dimension reconstruction approach, which makes it difficult to ensure compact reconstruction in each one-dimensional step. Moreover, the implementation becomes significantly complex in three-dimensional applications. 
Furthermore, in order to handle computations involving more complicated boundaries or flows requiring local mesh refinement, the development of high-order schemes on non-Cartesian structured meshes is also essential. Previous studies have developed high-order GKS in curvilinear coordinates using the dimension-by-dimension reconstruction method for both laminar \cite{HGKS-curv-1} and turbulent flows \cite{HGKS-curv-2}. However, due to the increased algorithmic complexity, the dimension-by-dimension approach still introduces significant challenges when extended to the compact scheme for three-dimensional cases.
Therefore, the development of a CGKS for three-dimensional generalized structured meshes is of critical importance, aiming to achieve high resolution and computational efficiency while preserving simplicity in implementation.

In this paper, we present a high-resolution compact gas-kinetic scheme designed for three-dimensional non-orthogonal meshes.
For this scheme, a novel multidimensional fifth-order compact reconstruction method is introduced for Cartesian meshes firstly. This method incorporates line-averaged derivatives within each control volume, offering additional degrees of freedom for constructing higher-order reconstruction by a compact stencil. It effectively enhances the resolution of CGKS.
Furthermore, to extend the new CGKS to non-orthogonal meshes, a simple coordinate transformation approach is employed. This method significantly reduces memory requirements and enhances computational efficiency compared to multidimensional reconstructions on the original meshes and dimension-by-dimension high-order methods.
Moreover, an efficient nonlinear combination method is employed to ensure that the compact scheme achieves both high resolution and numerical stability.
In addition, the proposed scheme is implemented on multiple graphics processing units (GPUs) using the Compute Unified Device Architecture (CUDA). 
Numerical tests are performed on 8 NVIDIA GeForce RTX 4090 GPUs to validate the scheme and accelerate computations.

This paper is organized as follows: Section 2 introduces the gas-kinetic scheme and finite volume framework. Section 3 presents the new fifth-order compact gas-kinetic scheme on Cartesian meshes, while Section 4 extends the compact scheme to non-orthogonal meshes and presents its parallel implementation. Numerical examples are provided in Section 5, and conclusions are drawn in the final section.

\section{Gas-kinetic scheme and finite volume method}

In the gas-kinetic method, a physically modeled, space-time coupled solution of the BGK equation for the gas distribution function is employed to compute time-dependent numerical fluxes \cite{GKS-Xu1,GKS-Xu2}.
The BGK equation is a simplification of Boltzmann equation, and the three-dimensional BGK equation \cite{BGK-1,BGK-2} can be written as
\begin{equation}\label{bgk}
f_t+\boldsymbol{u}\cdot \nabla f=\frac{g-f}{\tau},
\end{equation}
where $\boldsymbol{u}=(u,v,w)$ is the particle velocity, $\tau$ is the collision time and $f$ is the gas distribution function. The
equilibrium state $g$ is given by Maxwellian distribution 
\begin{equation}
g=\rho\left(\frac{\lambda}{\pi}\right)^{(N+3)/2}e^{-\lambda[(\boldsymbol{u}-\boldsymbol{U})^2+\boldsymbol{\xi}^2]},
\end{equation}
where $\rho$ is the density, $\boldsymbol{U}=(U,V,W)$ is the macroscopic fluid velocity, and $\lambda=\displaystyle\frac{\rho}{2p}$, $p$ is the pressure. 
In the BGK equation, the collision operator involves a simple relaxation from $f$ to the local equilibrium state $g$. 
The variable $\boldsymbol{\xi}$ accounts for the internal degree of freedom, $\boldsymbol{\xi}^2=\xi_1^2+\dots+\xi_N^2$, where $N=(5-3\gamma)/(\gamma-1)$ is the number of
internal degree of freedom and $\gamma$ is the specific heat ratio. 
The collision term satisfies the compatibility condition
\begin{equation*}
\int \frac{g-f}{\tau}\boldsymbol{\psi} \, \text{d}\Xi=0,
\end{equation*}
where $\displaystyle\boldsymbol{\psi}=\left(1,u,v,w,\frac{1}{2}(u^2+v^2+w^2+\boldsymbol{\xi}^2)\right)^T$ and $\text{d}\Xi=\text{d}u\text{d}v\text{d}w\text{d}\xi_1\dots\text{d}\xi_{N}$.
The macroscopic conservative variables $\boldsymbol{Q}=(\rho, \rho U,\rho V, \rho W, \rho E)$ can be calculated through the gas distribution function $f$
\begin{equation}\label{Q-macro}
\boldsymbol{Q}=\int f \boldsymbol{\psi} \,\text{d}\Xi,
\end{equation}
and the corresponding fluxes can be given by taking moments of the gas
distribution function
\begin{align}\label{flux-macro}
\boldsymbol{F}=\int \boldsymbol{u} f \boldsymbol{\psi} \, \text{d}\Xi,
\end{align}
where $\boldsymbol{F}$ represents either the Euler flux or the Navier-Stokes (NS) flux, depending on the relationship between $f$ and $g$ \cite{GKS-Xu1,GKS-Xu2}.
In GKS, $\boldsymbol{F}$ at cell interface is determined by the gas distribution function $f$.
Based on the integral solution of BGK equation 
Eq.\eqref{bgk}, the gas distribution function $f(\boldsymbol{x}_{G},t,\boldsymbol{u},\xi)$ can be given by
\begin{equation*}
f(\boldsymbol{x}_{G},t,\boldsymbol{u},\xi)=\frac{1}{\tau}\int_0^t
g(\boldsymbol{x}',t',\boldsymbol{u}, \xi)e^{-(t-t')/\tau}\text{d}t'+e^{-t/\tau}f_0(-\boldsymbol{u}t,\xi),
\end{equation*}
where $\boldsymbol{x}'=\boldsymbol{x}_{G}-\boldsymbol{u}(t-t')$ is the trajectory of particles and $f_0$ is the initial gas distribution function. 
In order to calculate the flux $\boldsymbol{F}$ by $f$ and used in the numerical schemes, the integral solution is modeled \cite{GKS-Xu2}. The second-order approximate solution of the integral solution is
\begin{align}\label{flux}
f(\boldsymbol{x}_{G},t,\boldsymbol{u},\xi)=&(1-e^{-t/\tau})g_0+
[(t+\tau)e^{-t/\tau}-\tau](\overline{a}_1u+\overline{a}_2v+\overline{a}_3w)g_0\nonumber\\
+&(t-\tau+\tau e^{-t/\tau}){\bar{A}} g_0\nonumber\\
+&e^{-t/\tau}g_r[1-(\tau+t)(a_{1}^{r}u+a_{2}^{r}v+a_{3}^{r}w)-\tau A^r](1-H(u))\nonumber\\
+&e^{-t/\tau}g_l[1-(\tau+t)(a_{1}^{l}u+a_{2}^{l}v+a_{3}^{l}w)-\tau A^l]H(u),
\end{align}
where the terms related to $g_{0}$ are from the integral of the equilibrium state and the terms related to $g_l$ and $g_r$ are from the initial term $f_0$.
More details of the gas-kinetic scheme can be found in \cite{GKS-Xu1,GKS-Xu2}. 

For the second-order evolution solution, the time accurate distribution function $f(t)$ at cell interface can be approximated through a linearization in time \cite{CGKS-high-5}
\begin{equation*}
\hat{f}(t)=f^n+t f_t^n.
\end{equation*}
The two coefficients $f^n$ and $f_t^n$ are calculated as follows 
\begin{align*}
    f^n&=\big(4\bar{f}(\Delta t/2) - \bar{f}(\Delta t)\big)/\Delta t,\\
    f_t^n&=4\big(\bar{f}(\Delta t) - 2\bar{f}(\Delta t/2)\big)/{\Delta t}^2,
\end{align*}
where $\bar{f}(\Delta t)$ and $\bar{f}(\Delta t/2)$ are the time integrations of $f(t)$ over the interval $[t_n, t_n + \Delta t]$ and $[t_n, t_n + \Delta t/2]$, respectively.
The numerical fluxes and their time derivatives can be obtained by taking moments of $\hat{f}(t)$ and $\hat{f}_t(t)$ at $t = t_n$
\begin{equation}\label{deri-t-f}
\boldsymbol{F}^n = \int \boldsymbol{u} f^n \boldsymbol{\psi} \, \text{d}\Xi,\quad
\boldsymbol{F}_t^n = \int \boldsymbol{u} f_t^n \boldsymbol{\psi} \, \text{d}\Xi.
\end{equation}
Simultaneously, the flow variables and their time derivatives can be given by
\begin{equation}\label{deri-t-q}
\boldsymbol{Q}^n = \int f^n \boldsymbol{\psi} \, \text{d}\Xi,\quad
\boldsymbol{Q}_t^n = \int f_t^n \boldsymbol{\psi} \, \text{d}\Xi.
\end{equation}
The above quantities will be utilized to update the flow variables.

Taking moments of Eq.\eqref{bgk} and integrating with respect to space, the semi-discretized finite volume scheme can be obtained as
\begin{align}\label{semi}
\frac{\text{d} \boldsymbol{Q}_i}{\text{d} t}=\mathcal{L}(\boldsymbol{Q}_{i}),
\end{align}
where $\boldsymbol{Q}_i$ is the cell-averaged conservative variables over cell $\Omega_{i}$. 
The operator $\mathcal{L}$ is defined as
\begin{equation*}%\label{finite}
\mathcal{L}(\boldsymbol{Q}_{i})%=-\frac{1}{|\Omega_{i}|}\sum_{i_p\in N(i)}F_{i,i_p}(t)S_{i_p}
=-\frac{1}{|\Omega_{i}|}\sum_{j=1}^{N_i} \int_{\Sigma_{j}}\boldsymbol{F}(t)\boldsymbol{n}\,\text{d}S,
\end{equation*}
where $|\Omega_{i}|$ is the volume of $\Omega_{i}$, $\Sigma_{j}$ is the cell interfaces of $\Omega_{i}$, $N_i$ is the total number of faces of $\Omega_{i}$, and $\boldsymbol{n}=( \boldsymbol{n_x}, \boldsymbol{n_y}, \boldsymbol{n_z})$ is the unit direction of the cell interface. 
To achieve the expected order of accuracy, the Gaussian 
quadrature is used for the flux integration
\begin{align*}
\int_{\Sigma_{j}}\boldsymbol{F}(t)\boldsymbol{n}\,\text{d}S=\sum_{j=1}^{N_i}\sum_{k=1}^{K}\omega_{k}\boldsymbol{F}(\boldsymbol{x}_k,t)S_{j},
\end{align*}
where $S_{j}$ is the area of $\Sigma_{j}$, $\omega_{k}$ is the Gaussian quadrature weight, $\boldsymbol{x}_{k}$ is the Gaussian quadrature point, $K$ is the total number of Gaussian quadrature points.
In this study, structured meshes are utilized, where each face is a quadrilateral with $K=4$ Gaussian points.
The numerical flux $\boldsymbol{F}(\boldsymbol{x}_k,t)$ at Gaussian quadrature point can be given by Eq.\eqref{flux-macro}.

\section{Compact gas-kinetic scheme on Cartesian meshes}

The time-dependent gas distribution function of GKS not only enables the computation of numerical fluxes but also provides pointwise values of conservative variables at the cell interfaces. Consequently, the CGKS can simultaneously update both the cell-averaged conservative variables and their corresponding cell-averaged derivatives within each time step, without requiring additional equations to compute the cell-averaged gradients. 
In previous studies, the application of CGKS to three-dimensional structured meshes remains challenging due to the inherent complexity of the algorithm using a dimension-by-dimension reconstruction approach.
Therefore, it is essential to develop a compact GKS for three-dimensional structured meshes that not only achieves high accuracy and efficiency but also provides high-resolution performance. In this section, a new multidimensional compact reconstruction method is proposed for a novel CGKS on Cartesian meshes.

\subsection{Flow variables update in CGKS}\label{variable}

For the compact GKS, the cell-averaged conservative variable $\boldsymbol{Q}(t)$ over a control volume $\Omega$ is updated through the conservation law.
Simultaneously, the cell-averaged derivatives can be updated based on the pointwise values of the conservative variables at the cell interfaces.
By applying Gauss's theorem, the cell-averaged gradient of the conservative variable $\boldsymbol{Q}(t)$ over the control volume $\Omega$ can be expressed as
\begin{equation}\label{derivative}
|\Omega|(\nabla \boldsymbol{Q})(t)=\int_{\Omega}\nabla \boldsymbol{Q}(t)\text{d}V=\int_{\partial\Omega}\boldsymbol{Q}(t) \boldsymbol{n}\,\text{d}S,
\end{equation}
where $\boldsymbol{n}=( \boldsymbol{n_x}, \boldsymbol{n_y}, \boldsymbol{n_z})$ is the unit direction on the cell interface. The variable $\boldsymbol{Q}$ at the cell interface is computed through Eq.~\eqref{Q-macro}, and the Gaussian quadrature is employed to evaluate the surface integral of $\boldsymbol{Q}(t)$ over the cell interface.

In this work, the line-averaged derivatives within each cell are introduced to construct a higher-order CGKS with enhanced resolution on a compact stencil. 
These line-averaged derivatives act as additional internal degrees of freedom, which further improve the overall resolution of the scheme.
The line-averaged partial derivative along a line segment defined by the points $\boldsymbol{x}_{1}$ and $\boldsymbol{x}_{2}$ is given by
\begin{equation}\label{derivative-l}
\begin{aligned}
(\partial_{l} \boldsymbol{Q})(t)&=\frac{1}{\left|l_G\right|}\int_{\boldsymbol{x}_{1}}^{\boldsymbol{x}_{2}}\frac{\partial \boldsymbol{Q}}{\partial l}(t)\text{d}l\\
&=\frac{1}{\left|l_G\right|}\left(\boldsymbol{Q}(\boldsymbol{x}_{2},t)-\boldsymbol{Q}(\boldsymbol{x}_{1},t)\right),
\end{aligned}
\end{equation}
where $\left|l_G\right|=\left|\boldsymbol{x}_{2}-\boldsymbol{x}_{1}\right|$, $\boldsymbol{Q}(\boldsymbol{x}_{1},t)$ and $\boldsymbol{Q}(\boldsymbol{x}_{2},t)$ are the conservative variables on this two points, respectively.
Thus, the cell-averaged conservative variable, the cell-averaged gradient and line-averaged derivatives are utilized to perform spatial reconstruction and are updated simultaneously.

\begin{figure}[!h]
\centering
\includegraphics[width=0.4\textwidth]{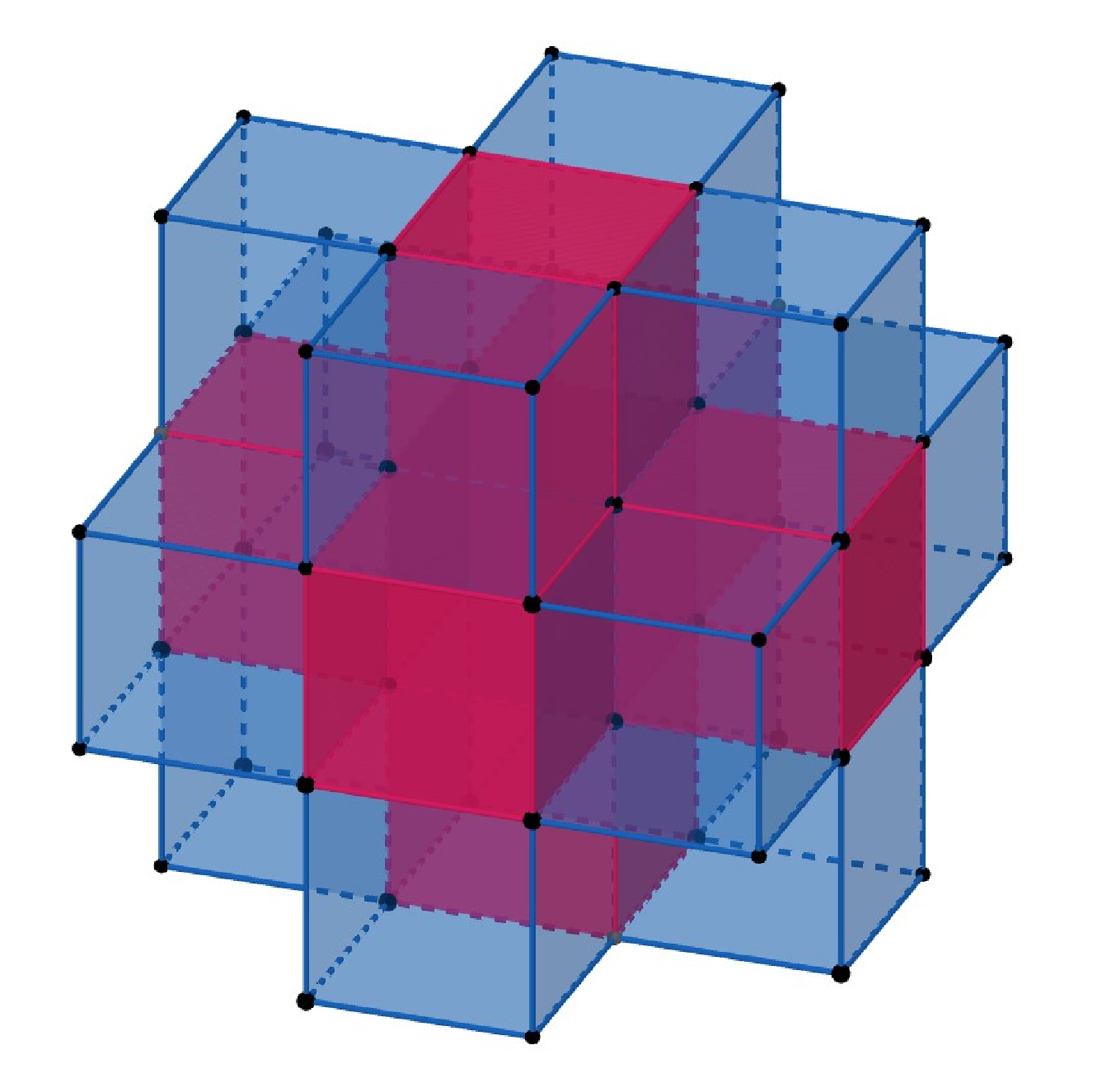}{(a)}
\includegraphics[width=0.445\textwidth]{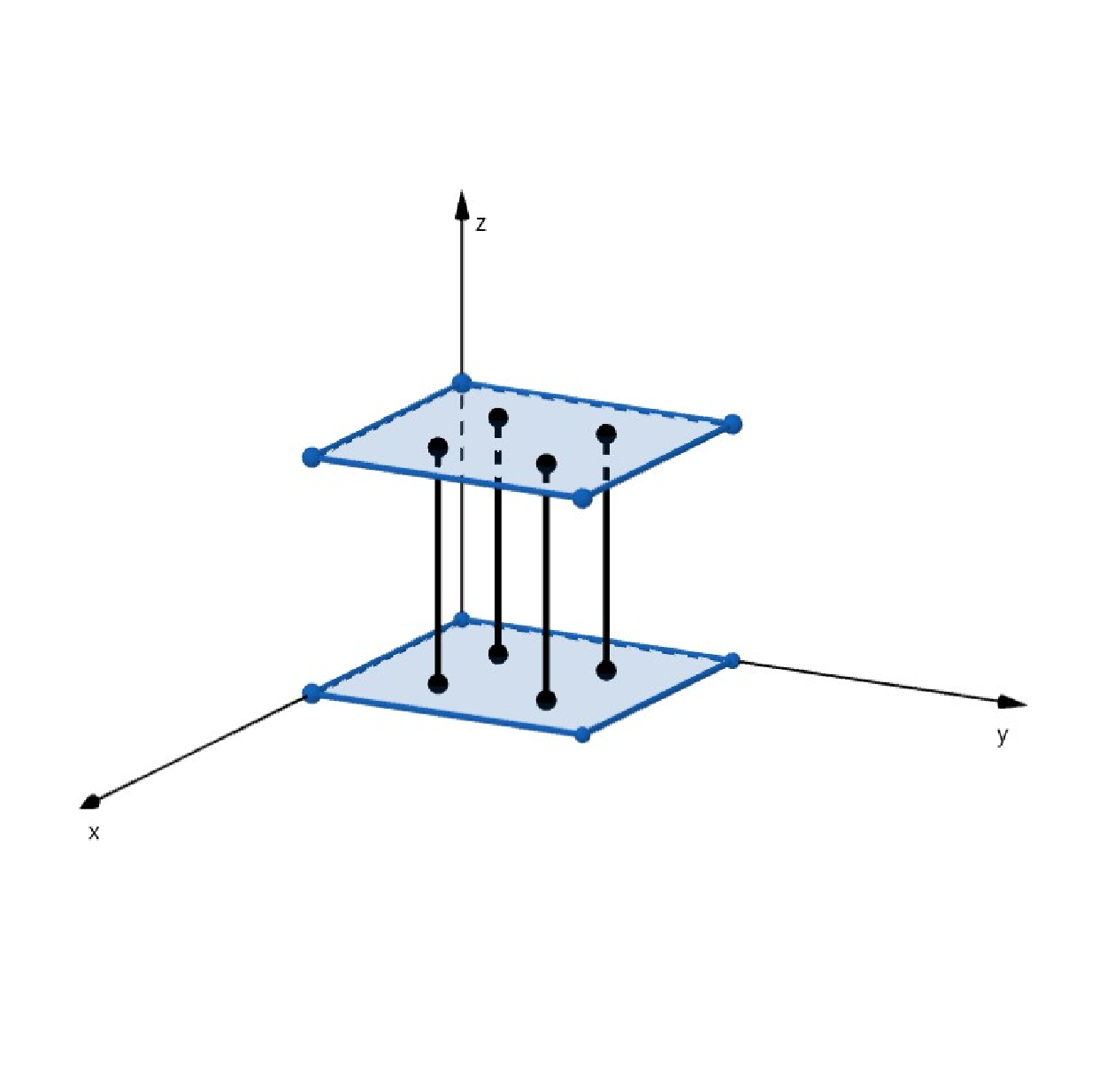}{(b)}
\caption{\label{stencil}(a) The schematic diagram of the compact stencil used for spatial reconstruction in the fifth-order CGKS. (b) The schematic diagram of the line-averaged partial derivatives in the z-direction within a cell.}
\end{figure}

\subsection{Stencils for the compact reconstruction}

For the target cell $\Omega_0 = \Omega_{ijk}$ in a three-dimensional Cartesian mesh, Figure \ref{stencil}(a) presents the schematic of the compact stencil employed for spatial reconstruction in the new fifth-order CGKS. The stencil consists of the face-neighboring cells (red cubes) and edge-neighboring cells (blue cubes) surrounding the target cell.
To facilitate the description, the following definitions are introduced. For the target cell $\Omega_0$, the cell interfaces are labeled as $F_m$, $m = 1, \dots, 6$, and the cell edges are labeled as $E_n$, $n = 1, \dots, 12$. The face-neighboring cell of $\Omega_0$ that shares the face $F_m$ is denoted by $\Omega_{1,m}$, while the edge-neighboring cell that shares only the edge $E_n$ is denoted by $\Omega_{2,n}$.

As shown in Figure \ref{stencil}(a), the big stencil used to construct the fifth-order reconstruction polynomial is defined as follows
\begin{equation*}
S^{cell-5}=\{\Omega_0, \Omega_{1,m}, \Omega_{2,n}\,|\,m = 1, \dots, 6; \, n = 1, \dots, 12\}.
\end{equation*}
As discussed in Section \ref{variable}, for each cell in the compact stencil $S^{cell-5}$, the cell-averaged variable, cell-averaged gradient, and line-averaged derivatives can be explicitly obtained. These quantities are then used to construct the fifth-order polynomial. This study proposes a novel reconstruction for the CGKS, designed to achieve the balance between compactness, resolution, and numerical stability. The scheme utilizes information from the target cell to enhance resolution, and its stability and robustness are validated through numerical experiments.
For achieving fifth-order accuracy, the reconstruction utilizes the following three types of data
\begin{align}\label{stencil-5}
S^{data-5}=S^{0}\cup S^{face}\cup S^{edge},
\end{align}
where
\begin{align*}
    S^{0}&=\{\boldsymbol{Q}_0,\,(\partial_l \boldsymbol{Q})_{0,G}\},\\
    S^{face}&=\{\boldsymbol{Q}_{1,m},\,(\nabla \boldsymbol{Q})_{1,m}\,|\,m = 1, \dots, 6\},\\
    S^{edge}&=\{\boldsymbol{Q}_{2,n},\,(\partial_{\boldsymbol{\tau}} \boldsymbol{Q})_{2,n}\,|\,n=1,\dots,12\}.
\end{align*}
Detailed descriptions of these quantities are provided below.

In $S^{data-5}$, $\boldsymbol{Q}_0$, $\boldsymbol{Q}_{1,m}$ and $\boldsymbol{Q}_{2,n}$ represent the cell-averaged conservative values over cells $\Omega_0$, $\Omega_{1,m}$ and $\Omega_{2,n}$, respectively. $(\nabla \boldsymbol{Q})_{1,m}$ denotes the cell-averaged gradient over cell $\Omega_{1,m}$, where $\nabla=(\partial_x,\partial_y,\partial_z)$. 
$(\partial_{\boldsymbol{\tau}} \boldsymbol{Q})_{2,n}$ denotes the cell-averaged directional derivative over cell $\Omega_{2,n}$, where $\partial_{\boldsymbol{\tau}}=(\partial_{\tau_1}, \partial_{\tau_2})$. Here, $\partial_{\tau_1}$ and $\partial_{\tau_2}$ are expressed as linear combinations of $\partial_x$, $\partial_y$ and $\partial_z$. The method of combination will be introduced in Section \ref{reconst-sec}.
$(\partial_{l} \boldsymbol{Q})_{0}$ represents the line-averaged partial derivative within the target cell $\Omega_0$. 
For the two opposite faces of the target cell along the same coordinate axis, four Gaussian points are selected on each face. Each Gaussian point on one face is uniquely paired with a corresponding Gaussian point on the opposite face, and these pairs are connected by straight lines. $(\partial_{l} \boldsymbol{Q})_{0,G}$ denotes the line-averaged partial derivative computed along one of these connecting lines.
For example, Figure \ref{stencil} (b) illustrates the schematic diagram of the line-averaged partial derivatives in the z-direction.
Consequently, there are a total of 12 line-averaged partial derivatives within the target cell. For convenience, all the line-averaged derivatives are rearranged as $\{(\partial_{l} \boldsymbol{Q})_{s},\,s=1,\dots,12\}$.

To address discontinuities, such as shock waves that appear in compressible flows, nonlinear reconstruction is often required. In this study, the nonlinear reconstruction is achieved through a weighted nonlinear combination of fifth-order and linear polynomials.
To construct several linear polynomials, sub-stencils $S_{m}^{cell-2}$ ($m = 1, \dots, 6$) for the target cell $\Omega_0$ are selected as follows
\begin{equation*}
S_{m}^{cell-2}=\{\Omega_0, \Omega_{1,m}\}.
\end{equation*}
For second-order accuracy, the reconstruction employs the following two types of data
\begin{align}\label{stencil-2}
S_{m}^{data-2}=S_{m}^{0}\cup S_{m}^{face},\,m=1,\dots,6,
\end{align}
where
\begin{align*}
S_{m}^{0}&=\{\boldsymbol{Q}_0,\,(\partial_l \boldsymbol{Q})^m_{0,G}\},\\
S_{m}^{face}&=\{\boldsymbol{Q}_{1,m}\},   %,\,(\nabla \boldsymbol{Q})_{1,m}
\end{align*}
$\boldsymbol{Q}_{1,m}$ is the cell-averaged conservative value over cell $\Omega_{1,m}$, $(\partial_l \boldsymbol{Q})^m_{0,G}$ is the selected line-averaged partial derivative within the target cell $\Omega_0$.
% For the sub-stencil $S_{m}^{CGKS-2}$, among the 12 line-averaged positions of the target cell defined earlier, 8 line-averaged positions are selected. 
% These positions are chosen such that the lines are not parallel to the vector connecting the centroids of the neighboring cell $\Omega_{1,m}$ and the target cell $\Omega_0$.
Eight out of the total twelve line-averaged derivatives are utilized for the reconstruction of linear polynomials. The directions of these eight directional derivatives are perpendicular to the connecting line between the centroids of two adjacent cells in the sub-stencil. 
%This choice enhances the stability of the scheme. % and improves its robustness in capturing shock waves
For convenience, all the line-averaged derivatives are rearranged as $\{(\partial_{l} \boldsymbol{Q})_{s^{'}},\,s^{'}=1,\dots,8\}$.
Compared to traditional reconstruction methods, the stencils presented in this section increase the degrees of freedom by incorporating line-averaged gradients of the target cell. This approach utilizes more detailed information to capture gradient variations, thereby further improving resolution.

\subsection{Compact reconstruction}\label{reconst-sec}

For the target cell $\Omega_0$, a quartic polynomial $P^4(\boldsymbol{x})$ can be constructed based on the big stencil $S^{CGKS-5}$ in Eq.\eqref{stencil-5}.
The quartic polynomial $P^4(\boldsymbol{x})$ is defined as
\begin{equation*}
P^4(\boldsymbol{x})=\boldsymbol{Q}_{0}+\sum_{|\boldsymbol{d}|=1}^4 a_{\boldsymbol{d}}p_{\boldsymbol{d}}(\boldsymbol{x}),
\end{equation*}
where $\boldsymbol{Q}_{0}$ is the cell-averaged variable over cell $\Omega_{0}$, the multi-index $\boldsymbol{d}=(d_1, d_2, d_3)$, and $|\boldsymbol{d}|=d_1+d_2+d_3$. 
The zero-mean basis $p_{\boldsymbol{d}}(\boldsymbol{x})$ is defined as
\begin{equation}\label{base}
p_{\boldsymbol d}(\boldsymbol{x})=\displaystyle\frac{1}{d_1!d_2!d_3!}{\delta x}^{d_1}{\delta y}^{d_2}{\delta z}^{d_3}-\frac{1}{\left|\Omega_{0}\right|}\int_{\Omega_{0}}\frac{1}{d_1!d_2!d_3!}{\delta x}^{d_1}{\delta y}^{d_2}{\delta z}^{d_3}\text{d}V,
\end{equation}
where
\begin{equation*}
\delta x=\frac{1}{h_x}(x-x_0),\,\delta y=\frac{1}{h_y}(y-y_0),\delta z=\frac{1}{h_z}(z-z_0),
\end{equation*}
$h_x$, $h_y$ and $h_z$ are the characteristic scales of $\Omega_0$ along the three directions of axes, and $(x_0,y_0,z_0)$ is the centroid of $\Omega_0$.
To determine the fifth-order polynomial, the following constrains need to be
satisfied
\begin{equation}\label{compact-big-5}
\begin{split}
\frac{1}{\left|\Omega_{k}\right|}\int_{\Omega_{k}}P^4(\boldsymbol{x})\text{d}V&=\boldsymbol{Q}_{k},~\Omega_{k}\in S^{cell-5},\\
\frac{h_1}{\left|\Omega_{1,m}\right|}\int_{\Omega_{1,m}}\nabla P^4(\boldsymbol{x})\text{d}V&=h_1\cdot(\nabla \boldsymbol{Q})_{1,m},~m=1,\dots, 6,\\
\frac{h_2}{\left|\Omega_{2,n}\right|}\int_{\Omega_{2,n}}\frac{\partial}{\partial_{\boldsymbol{\tau}}}P^4(\boldsymbol{x})\text{d}V&=h_2\cdot(\partial_{\boldsymbol{\tau}} \boldsymbol{Q})_{2,n},~n=1,\dots,12,\\
\frac{h_3}{\left|l_G\right|}\int_{l_G}\frac{\partial}{\partial_{l}}P^4(\boldsymbol{x})\text{d}l&=h_3\cdot(\partial_{l} \boldsymbol{Q})_{s},~s=1,\dots,12,
\end{split}
\end{equation}
where the parameters $h_1$, $h_2$, and $h_3$ are used for scaling, the details of which will be introduced later, $\boldsymbol{Q}_{k}$ is the cell-averaged variable over the cell $\Omega_{k}$, where $\Omega_{k} \in S^{cell-5}$. The cell-averaged gradient $(\nabla \boldsymbol{Q})_{1,m}$ over the cell $\Omega_{1,m}$ and the line-averaged derivative $(\partial_{l} \boldsymbol{Q})_{s}$ on the right-hand side are calculated as described in Section \ref{variable}.
For the edge-neighboring cell $\Omega_{2,n}$, which shares the common edge $E_{n}$, the operator $\partial_{\boldsymbol{\tau}} = (\partial_{\tau_1}, \partial_{\tau_2})$ corresponds to the directional derivatives along two distinct directions. One direction is parallel to the edge $E_{n}$, while the other is parallel to the vector connecting the centroid of the target cell to the centroid of $\Omega_{2,n}$. Consequently, the cell-averaged directional derivatives $(\partial_{\boldsymbol{\tau}} \boldsymbol{Q})_{2,n}$ can be computed as a linear combination of the derivatives $(\partial_x \boldsymbol{Q})_{2,n}$, $(\partial_y \boldsymbol{Q})_{2,n}$, and $(\partial_z \boldsymbol{Q})_{2,n}$.
Through this selection process, the amount of derivative information required from the edge-neighboring cells is reduced, while still maintaining the over-determined solvability of the system described in Eq.~\eqref{compact-big-5}.

In Eq.~\eqref{compact-big-5}, the derivative terms are scaled by the parameters $h_1$, $h_2$, and $h_3$. These parameters serve to reduce the condition number of the reconstruction matrix and improve the stability and resolution of the final compact scheme.
In previous constructions of compact schemes, $h_1$, $h_2$, and $h_3$ were all set to $h$ \cite{CGKS-high-4}, where $h$ is the characteristic scale of the cell $\Omega_0$. 
This choice ensured that the elements in each row of the same column in the reconstruction matrix had the same order of magnitude, thereby optimizing the condition number of the matrix. 
In this study, more refined adjustments have been made to $h_1$, $h_2$, and $h_3$ to enhance resolution. Specifically, we set $h_1 = 2h$, $h_2 = h$, and $h_3 = 3h$. The parameters $h_3$, $h_1$ and $h_2$ decrease sequentially, reflecting the progressively diminishing importance of the relationships between these three types of data and the target cell.
% In detail, the line-averaged derivatives from $S^{0}$ are assigned the highest weight because they represent the internal degrees of freedom of the target cell, which play a crucial role in achieving higher resolution. The gradients provided by $S^{face}$ are derived from the neighboring cells closest to the target cell, so they are assigned the second-highest weight to further promote compactness. In contrast, the derivatives from $S^{edge}$ are assigned the lowest weight, ensuring that the system of equations remains overdetermined.
Furthermore, in order to ensure linear stability, the system in Eq.\eqref{compact-big-5} is solved using the constrained least-square method. Only the conservative variable equations for $\{\boldsymbol{Q}_{0}, \boldsymbol{Q}_{1,m} | m=1,\dots,6\}$ are set as strictly satisfied and others are satisfied in the least-square sense.

For the target cell $\Omega_0$, several linear polynomials $P_m^1(\boldsymbol{x})$ can be constructed based on $S_{m}^{data-2}$ ($m=1,\dots,6$) in Eq.\eqref{stencil-2}.
The linear polynomial $P_m^1(\boldsymbol{x})$ is defined as
\begin{equation*}
P_m(\boldsymbol{x})=\boldsymbol{Q}_{0}+\sum_{|\boldsymbol d|=1}b_{\boldsymbol d}^mp_{\boldsymbol d}(\boldsymbol{x}),
\end{equation*}
where $m=1,\dots,6$, and the basis function $p_{\boldsymbol d}(\boldsymbol{x})$ is defined as same as Eq.\eqref{base}.
To determine these linear polynomials, the following constrains need to be
satisfied
\begin{equation}\label{compact-sub-2}
\begin{split}
\frac{1}{\left|\Omega_{1,m}\right|}\int_{\Omega_{1,m}}P_m^1(\boldsymbol{x})\text{d}V&=\boldsymbol{Q}_{1,m},\\
\frac{1}{\left|l_G\right|}\int_{l_G}\frac{\partial}{\partial_{l}}P_m^1(\boldsymbol{x})\text{d}l&=(\partial_{l} \boldsymbol{Q})_{s^{'}},~s^{'}=1,\dots,8,
\end{split}
\end{equation}
where $\boldsymbol{Q}_{1,m}$ is the cell-averaged variable over cell $\Omega_{1,m}$, and line-averaged partial derivative $(\partial_{l} \boldsymbol{Q})_{s^{'}}$ is calculated in the same way as presented in Eq.\eqref{derivative-l}.
The constrained least-square method is used for solving the above linear systems, where the conservative variable equations are set as strictly satisfied and others are satisfied in the least-square sense.

%================================================================================================
%================================================================================================

To solve problems involving discontinuities, an efficient nonlinear combination of high-order and low-order polynomials is employed. 
In this study, the generalized ENO (GENO) nonlinear reconstruction method \cite{GENO} is adopted. The GENO method adaptively transitions from high-order linear reconstruction in smooth flow regions to second-order reconstruction near discontinuities, achieving an optimal balance between high accuracy (maintained through high-order linear reconstruction across a wide wavenumber range) and shock-capturing robustness (ensured by preserving the ENO property). 
With the reconstructed polynomial $P^4(\boldsymbol{x})$ and $P_m^1(\boldsymbol{x}), m=1,...,6$, the point-value $\boldsymbol{Q}(\boldsymbol{x}_{G})$ and the spatial derivatives
$\partial_{x,y,z} \boldsymbol{Q}(\boldsymbol{x}_{G})$ at Gaussian quadrature point for the nonlinear reconstruction are provided as follows.
\begin{equation}\label{weno-new}
\begin{split}
\boldsymbol{Q}(\boldsymbol{x}_{G})&=\chi P^4(\boldsymbol{x}_{G}) + (1-\chi)\left(\sum_{m=1}^{6}\omega_{m} P_m^1(\boldsymbol{x})\right),\\
\partial_{x,y,z} \boldsymbol{Q}(\boldsymbol{x}_{G})&=\chi \partial_{x,y,z} P^4(\boldsymbol{x}_{G}) + (1-\chi)\left(\sum_{m=1}^{6}\omega_{m} \partial_{x,y,z} P_m^1(\boldsymbol{x})\right),
\end{split}
\end{equation}
where 
\begin{equation*}
\displaystyle\omega_{m}=\frac{\overline{\omega_{m}}}{\sum_{m=1}^{6} \overline{\omega_{m}}},\,
\displaystyle\overline{\omega_{m}}=\frac{d_m}{(IS_m+\epsilon)^5},
\end{equation*}
and $d_1=\dots=d_6=\displaystyle\frac{1}{6}$ are the linear weights, $\epsilon$ is a small positive constant, set to $10^{-15}$.
In the present study, a minor modification is introduced to the path function, $\chi$, of the GENO reconstruction scheme \cite{GENO}. The path function $\chi$ is originally determined by the smoothness indicators $IS^H$ and $IS^L$, which correspond to the high-order linear polynomial and smooth polynomials from the low-order ENO reconstruction, respectively \cite{GENO}.
While, in the present study, $IS^L$ is also given by the high-order linear reconstruction polynomial as
\begin{equation*}
    IS^L=\sum_{|l|=1}h^{2|l|-3}\int_{\Omega}\big(D^{|l|}P^4(\boldsymbol{x})\big)^2\,\mathrm{d}\Omega.
\end{equation*}
Notably, this calculation for $IS^L$ exclusively involves terms related to the first-order derivatives. Consequently, the path function becomes
\begin{equation*}
    \chi=\tanh(C\alpha)/\tanh(C),
\end{equation*}
where 
\begin{equation*}
    \begin{aligned}
    \alpha&=\displaystyle\frac{2\alpha^H}{\alpha^H+\alpha^L},\\
    \alpha^H&=\displaystyle 1+(\frac{IS^{\tau}}{IS^H+\epsilon})^r,\,\alpha^L=\displaystyle 1+(\frac{IS^{\tau}}{IS^L+\epsilon})^r,
    \end{aligned}
\end{equation*}
and the parameter $r$ is set to 2. The computation of other variables in the GENO reconstruction of this study remains consistent with the multi-dimensional reconstruction formulas presented in \cite{GENO}.

\begin{figure}[!h]
\centering
\includegraphics[width=0.4\textwidth]{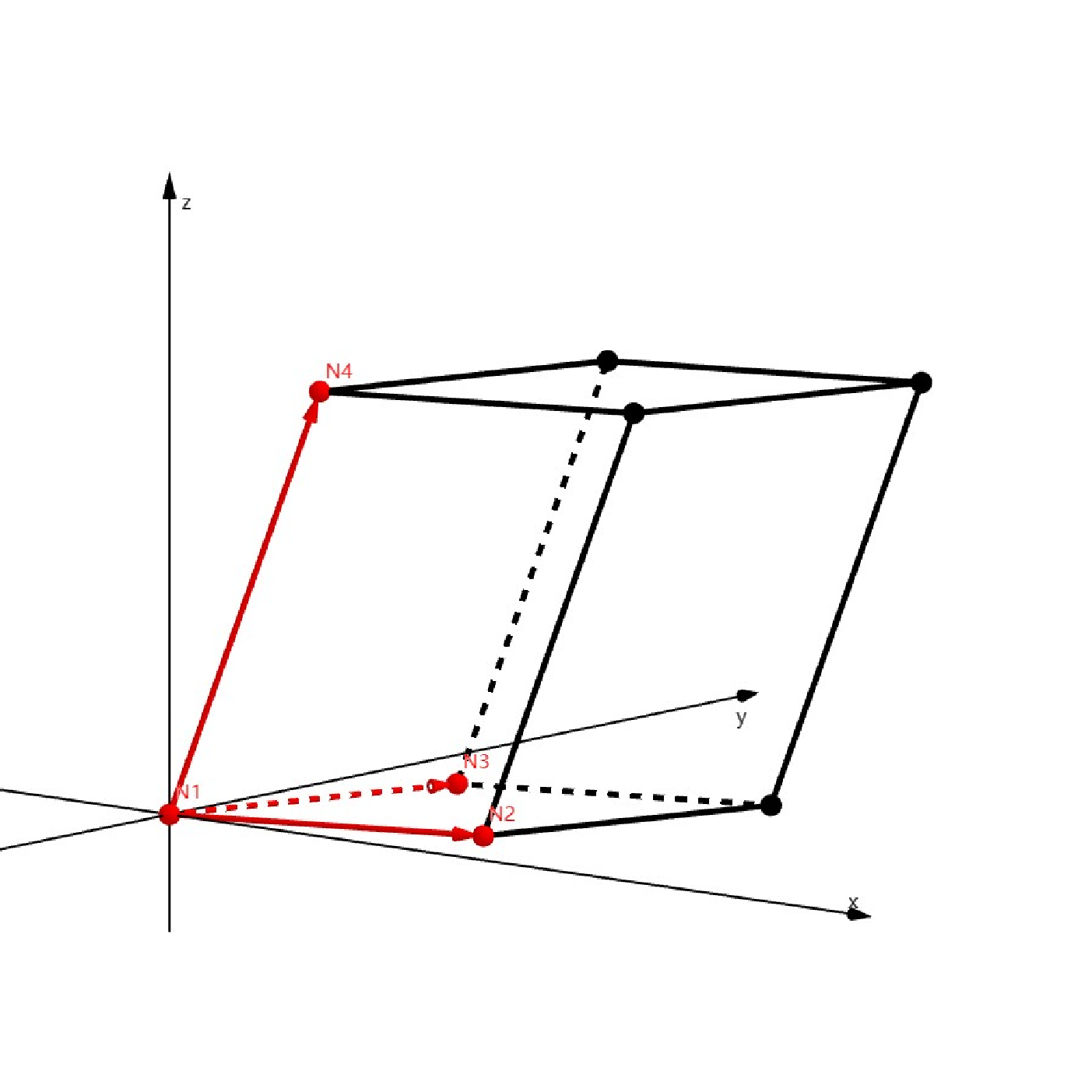}{(a)}
\includegraphics[width=0.4\textwidth]{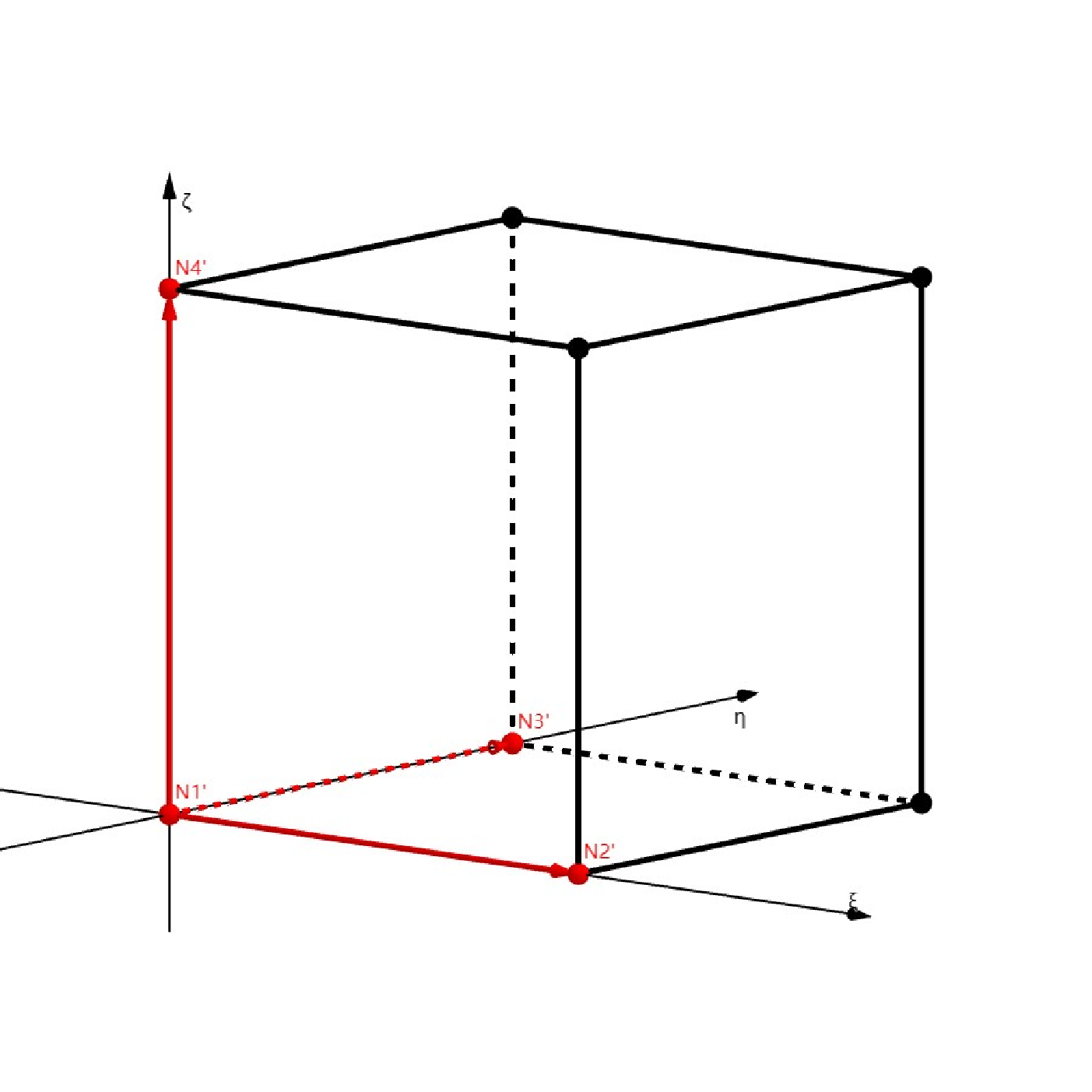}{(b)}
\caption{\label{trans} (a) The target cell $\Omega_0$ in the physical coordinate system $(x, y, z)$. (b) The standard cell ${\Omega_0}^{'}$ in the reference coordinate system $(\xi, \eta, \zeta)$.}
\end{figure}

\section{Compact gas-kinetic scheme on non-orthogonal meshes}

To address computations involving curved boundaries or flows that require local mesh refinement, the new fifth-order reconstruction must be extended to non-orthogonal meshes.
In this section, a straightforward transformation approach is introduced, offering advantages such as low memory requirements and simplicity, making it particularly suitable for large-scale computations.

For the target cell $\Omega_0$ in the physical coordinate system $(x, y, z)$, three vectors formed by the four nodes of $\Omega_0$ are selected, as shown by the red vectors in Figure \ref{trans}(a). Let the coordinates of these four nodes be $N_k=(x_k,y_k,z_k)$, $k=1,\cdots,4$. 
The process of transforming $\Omega_0$ into a standard cell ${\Omega_0}^{'}$ is achieved by converting these three vectors into the basis vectors of the reference coordinate system $(\xi, \eta, \zeta)$. As shown in Figure \ref{trans}(b), the coordinates of $N_k'$, $k=1,\cdots,4$, are $(0,0,0)$, $(1,0,0)$, $(0,1,0)$ and $(0,0,1)$, respectively.
The above transformation from the physical coordinate system $(x, y, z)$ into the reference coordinate system $(\xi, \eta, \zeta)$ is defined as
\begin{equation}
\left(
    \begin{array}{c}
    x\\
    y\\
    z\\
    \end{array}
\right)
=
\left(
    \begin{array}{c}
    x_1\\
    y_1\\
    z_1\\
    \end{array}
\right)
+
\boldsymbol{J}_0
\left(
    \begin{array}{c}
    \xi\\
    \eta\\
    \zeta\\
    \end{array}
\right),
\end{equation}
where $\boldsymbol{J}_0$ is the Jacobian matrix, and $\boldsymbol{J}_0$ takes
\begin{equation}
\boldsymbol{J}_0=
\left(
    \begin{array}{ccc}
    x_2-x_1 & x_3-x_1 & x_4-x_1\\
    y_2-y_1 & y_3-y_1 & y_4-y_1\\
    z_2-z_1 & z_3-z_1 & z_4-z_1\\
    \end{array}
\right).
\end{equation}
Based on the above Jacobian matrix $\boldsymbol{J}_0$, the face-neighboring cells $\Omega_{1,m}$ ($m=1,\dots,6$) and edge-neighboring cells $\Omega_{2,n}$ ($n=1,\dots,12$) are transformed to ${\Omega_{1,m}}^{'}$ and ${\Omega_{2,n}}^{'}$, respectively.

After the transformation from the target cell $\Omega_0$ to the standard cell ${\Omega_0}^{'}$, the reconstruction polynomial $P^4(\boldsymbol{\xi})$ for cell ${\Omega_0}^{'}$ can be constructed on the transformed cells 
\begin{equation*}
    {S^{cell-5}}^{'}=\{\Omega_0^{'}, \Omega_{1,m}^{'}, \Omega_{2,n}^{'}\,|\,m = 1, \dots, 6,\, n = 1, \dots, 12\}.
\end{equation*}
The quartic polynomial $P^4(\boldsymbol{\xi})$ is defined as
\begin{equation*}
P^4(\boldsymbol{\xi})=\boldsymbol{Q}_{0}^{'}+\sum_{|\boldsymbol{d}|=1}^4 a_{\boldsymbol{d}}^{'}p_{\boldsymbol{d}}(\boldsymbol{\xi}),
\end{equation*}
where $p_{\boldsymbol{d}}(\boldsymbol{\xi})$ is the zero-mean basis function defined on cell ${\Omega_0}^{'}$, and $\boldsymbol{Q}_{0}^{'}$ is the cell-averaged variable over cell ${\Omega_0}^{'}$.
To determine the polynomial $P^4(\boldsymbol{\xi})$, the system on transformed cells is written as
\begin{equation}\label{compact-big-ref}
\begin{split}
\frac{1}{\left|\Omega_{k}^{'}\right|}\int_{\Omega_{k}^{'}}P^4(\boldsymbol{\xi})\text{d}V^{'}&=\frac{\boldsymbol{Q}_k^{'}}{C_k},\,
~\Omega_k^{'}\in {S^{cell-5}}^{'},\\
\frac{h_1^{'}}{\left|\Omega_{1,m}^{'}\right|}\int_{\Omega_{1,m}^{'}}\nabla P^4(\boldsymbol{\xi})\text{d}V^{'}&=\frac{h_1^{'}}{C_{1,m}}\cdot(\nabla \boldsymbol{Q})_{1,m}^{'},\,
~m=1,\dots,6,\\
\frac{h_2^{'}}{\left|\Omega_{2,n}^{'}\right|}\int_{\Omega_{2,n}^{'}}\frac{\partial}{\partial_{\boldsymbol{\tau}}}P^4(\boldsymbol{\xi})\text{d}V^{'}&=\frac{h_2^{'}}{C_{2,n}}\cdot(\partial_{\boldsymbol{\tau}} \boldsymbol{Q})_{2,n}^{'},\,
~n=1,\dots,12,\\
h_3^{'}\cdot \int_{l_G^{'}}\frac{\partial}{\partial_{l}}P^4(\boldsymbol{\xi})\text{d}l^{'}&=h_3^{'}\cdot\left(\boldsymbol{Q}^{'}(\boldsymbol{x}_{G(s,2)})-\boldsymbol{Q}^{'}(\boldsymbol{x}_{G(s,1)})\right),\,~s=1,\dots,12,
%\left(P^4(\boldsymbol{\xi}_{G(s,2)})-P^4(\boldsymbol{\xi}_{G(s,1)})\right)
\end{split}
\end{equation}
where
\begin{equation*}
C_k=\frac{\left|\Omega_{k}^{'}\right|}{\left|\Omega_{k}\right|}|\boldsymbol{J}_0|,\quad \Omega_k^{'}\in {S^{cell-5}}^{'},
\end{equation*}
$\left|\Omega_{k}^{'}\right|$ is the volume of transformed cell $\Omega_{k}^{'}$, $\boldsymbol{J}_0$ is the Jacobian matrix of the target cell $\Omega_0$, $|\boldsymbol{J}_0|$ is the determinant of $\boldsymbol{J}_0$. Noted that the volume $\left|\Omega_{k}^{'}\right|$ of the transform cells should be calculated carefully. More details for the calculation can be found in the Appendix.
$\boldsymbol{Q}_k^{'}$ is the cell-averaged values over transformed cell $\Omega_{k}^{'}\in {S^{cell-5}}^{'}$, and $\boldsymbol{Q}_k^{'}=\boldsymbol{Q}_k$. 
$(\nabla \boldsymbol{Q})_{1,m}^{'}$ represents the cell-averaged gradients over cell ${\Omega_{1,m}}^{'}$ calculated by \cite{CGKS-high-5}
\begin{equation*}
(\nabla \boldsymbol{Q})_{1,m}^{'}(\xi,\eta,\zeta)=\boldsymbol{J}_0^T\cdot (\nabla \boldsymbol{Q})_{1,m}(x,y,z),
\end{equation*}
and $(\partial_{\boldsymbol{\tau}} \boldsymbol{Q})_{2,n}^{'}$ can be calculated in a similar way.
$\boldsymbol{x}_{G(s,1)}$ and $\boldsymbol{x}_{G(s,2)}$ represent a pair of one-to-one corresponding Gaussian points on the cell interfaces.
Due to the characteristic scale of cell $\Omega_0^{'}$ being 1, the scaling parameters are $h_1^{'}=2$, $h_2^{'}=1$, and $h_3^{'}=3$.
Similarly, the linear polynomials $P_m^1(\boldsymbol{\xi})$ ($m=1,\cdots,6$) for cell ${\Omega_0}^{'}$ can be constructed on the transformed cells.

With the reconstructed polynomial $P^4(\boldsymbol{\xi})$ and $P_m^1(\boldsymbol{\xi})$ ($m=1,\cdots,6$) the reconstructed point-value $\boldsymbol{Q}(\boldsymbol{\xi}_{G})^{'}$ and the spatial derivatives $\partial_{\xi,\eta,\zeta} \boldsymbol{Q}(\boldsymbol{\xi}_{G})^{'}$ in the reference coordinate system can be given by the nonlinear combination same as Eq.\eqref{weno-new}.
Then, the numerical flux at the interface is calculated after the reconstructed values $\boldsymbol{Q}(\boldsymbol{\xi}_{G})^{'}$ and $\partial_{\xi,\eta,\zeta} \boldsymbol{Q}(\boldsymbol{\xi}_{G})^{'}$ are transformed from the reference coordinate system back to the physical coordinate system.
The reverse transformation is accomplished by
\begin{align*}
\boldsymbol{Q}&=\boldsymbol{Q}^{'},\\
(\nabla \boldsymbol{Q})(x,y,z)&=(\boldsymbol{J}_0^T)^{-1}\cdot (\nabla \boldsymbol{Q})^{'}(\xi,\eta,\zeta).
\end{align*}

For the non-orthogonal meshes, multidimensional schemes typically require storing the coefficients $a_{\boldsymbol{d}}$ and $b_{\boldsymbol{d}}^m$ of the reconstruction polynomials for each computational cell.
The array that stores these coefficients has dimensions of $$O(N,\, \text{number of stencils},\, \text{number of basis functions},\, \text{amount of information}),$$ where $N$ is the total number of computational cells. In compact reconstruction, the amount of information is typically 4, which corresponds to the conservative variables and their gradients.
In large-scale computations, the storage requirements at this scale can become prohibitively expensive.
In this work, the proposed scheme utilizes a coordinate transformation that maps all target cells to a single standard cell. Consequently, only the coefficients of the standard cell need to be stored. 
This approach requires a matrix with dimensions $$O(1,\, \text{number of stencils}, \, \text{number of basis functions}, \, \text{amount of information}),$$
where the amount of information is $16$, representing the conservative variable, the cell-averaged gradient, and the line-averaged derivatives. 
Thus, the memory requirement is significantly reduced.
The above characteristic makes the current scheme ideal for large-scale parallel simulations, where reduced memory usage minimizes data access and transfer, thereby improving computational efficiency. Additionally, it enables larger simulations to be performed within limited memory resources.

\section{Temporal discretization}
    
In this section, the updates of the cell-averaged conservative variable, the cell-averaged gradient, and the line-averaged derivatives at each time step are briefly introduced.

Based on Eq.\eqref{semi}, the discrete conservation law in cell $\Omega_{i}$ is expressed as
\begin{equation*}
\boldsymbol{Q}_{i}^{n+1}=\boldsymbol{Q}_{i}^{n}+\int_{t^n}^{t^{n+1}} \mathcal{L}(\boldsymbol{Q}_{i},t)\,\mathrm{d}t,
\end{equation*}
where $\boldsymbol{Q}_{i}^{n+1}$ represents the cell-averaged conservative variable over cell $\Omega_{i}$ at $t_{n+1}=t_n+\Delta t$.
To achieve high-order temporal accuracy, the two-stage fourth-order temporal discretization method \cite{s2o4-0,GRP-high-1,GKS-high-1} is employed to update the cell-averaged conservative variable as follows
\begin{align*}
\boldsymbol{Q}_{i}^{*}&=\boldsymbol{Q}_{i}^n+\frac{1}{2}\Delta t\mathcal{L}(\boldsymbol{Q}_{i}^n)+\frac{1}{8}\Delta t^2\frac{\partial}{\partial t}\mathcal{L}(\boldsymbol{Q}^n_{i}), \\
\boldsymbol{Q}_{i}^{n+1}&=\boldsymbol{Q}_{i}^n+\Delta t\mathcal {L}(\boldsymbol{Q}_{i}^n)+\frac{1}{6}\Delta t^2\big(\frac{\partial}{\partial t}\mathcal{L}(\boldsymbol{Q}^n_{i})+2\frac{\partial}{\partial t}\mathcal{L}(\boldsymbol{Q}^{*}_{i})\big).
\end{align*}
As discussed in Section \ref{variable}, the cell-averaged gradient and the line-averaged derivatives are evaluated using the conservative variables at cell interfaces. The conservative variable at one side of a cell interface is evolved in two stages:
\begin{align*}
\boldsymbol{Q}^{*}&=\boldsymbol{Q}^n+\frac{1}{2}\Delta t (\partial_t \boldsymbol{Q})^n,\\
\boldsymbol{Q}^{n+1}&=\boldsymbol{Q}^n+\Delta t (\partial_t \boldsymbol{Q})^{*}.
\end{align*}
To provide variables on both sides of a cell interface, which may differ near discontinuities, the update model for $\boldsymbol{Q}(t)$ is expressed as
\begin{align*}
\boldsymbol{Q}^l(t)&=(1-e^{-\Delta t /\tau_0})\boldsymbol{Q}^e(t)+e^{-\Delta t /\tau_0}\boldsymbol{Q}_0^l(t),\\
\boldsymbol{Q}^r(t)&=(1-e^{-\Delta t /\tau_0})\boldsymbol{Q}^e(t)+e^{-\Delta t /\tau_0}\boldsymbol{Q}_0^r(t).
\end{align*}
Further calculation details can be found in \cite{CGKS-high-5,CGKS-high-4}.

\section{Numerical tests}

In this section, numerical tests ranging from subsonic to supersonic flows are conducted to evaluate the resolution, stability, and robustness of the proposed scheme.
For convenience, the fifth-order compact GKS is referred to as CGKS-5.
Unless otherwise specified, CGKS-5 indicates the scheme employing nonlinear reconstruction.

For evolution of flow fields, the time step $\Delta t$ is given by the CFL condition. In all test cases, CFL number is taken as $\text{CFL}=0.5$. 
For viscous flows, the time step is also limited by the viscous term as $\Delta t = \text{CFL}\cdot h^2/(3\nu)$, where $h$ is the cell size and $\nu$
is the kinematic viscosity coefficient. 
For the inviscid flows, the collision time $\tau$ takes
\begin{align*}
\tau=\epsilon \Delta t+C\displaystyle|\frac{p_l-p_r}{p_l+p_r}|\Delta t,
\end{align*}
where $p_l$ and $p_r$ denote the pressure on the left and right sides of the cell interface, $\epsilon=0.05$ and $C=5.0$. For the viscous flows, the collision time is related to the viscosity coefficient,
\begin{align*}
\tau=\frac{\mu}{p}+C \displaystyle|\frac{p_l-p_r}{p_l+p_r}|\Delta t,
\end{align*}
where $\mu$ is the dynamic viscous coefficient and $p$ is the pressure at the cell interface. In smooth flow regions, the collision time reduces to
\begin{equation*}
\tau=\frac{\mu}{p}.
\end{equation*}
In this section, the specific heat ratio $\gamma$ takes 1.4.

% Graphics Processing Units (GPUs), originally designed for graphics processing, have evolved into powerful tools for hardware acceleration, particularly well-suited for highly parallel computing tasks.
% Recently, GPU acceleration has been successfully applied to the high-order GKS, significantly enhancing its capability to tackle complex and large-scale computational problems \cite{GKS-para-2,GKS-para-3}.
For the numerical cases presented in this section, the CGKS-5 is implemented on multiple GPUs using the Compute Unified Device Architecture (CUDA) for parallel computation and the Message Passing Interface (MPI) for inter-process communication. 
This setup enables large-scale simulations, utilizing 8 NVIDIA GeForce RTX 4090 GPUs with double-precision.
The computational domain is partitioned into multiple subdomains along one dimension, with each subdomain assigned to a specific MPI process, which in turn is associated with a GPU. 
Each MPI process manages the parallel computation within its assigned subdomain by invoking CUDA kernels. 
% To maintain consistency across subdomains, boundary data is exchanged between processes through MPI inter-process communication functions.
In this work, the CGKS-5 demonstrates high parallelizability, efficiently executing kernels for tasks such as initialization, time step calculation, compact reconstruction, flux computation, and flow variable updating.

\begin{table}[!h]
\begin{center}
\def\temptablewidth{0.8\textwidth}{\rule{\temptablewidth}{1.0pt}}
\begin{tabular*}{\temptablewidth}{@{\extracolsep{\fill}}c|c|c|c|c} 
MESH            & $L^1$ Error & Order   & $L^2$ Error & Order  \\
\hline
$10^3$ & $2.1963\times 10^{-2}$ &  & $8.5136\times 10^{-3}$  & \\
$20^3$ & $7.1472\times 10^{-4}$ & 4.9415 & $2.8195\times 10^{-4}$  & 4.9163 \\
$40^3$ & $2.3127\times 10^{-5}$ & 4.9497 & $9.0788\times 10^{-6}$  & 4.9568 \\
$80^3$ & $8.0092\times 10^{-7}$ & 4.8518 & $3.1656\times 10^{-7}$  & 4.8420 \\
\end{tabular*}
{\rule{\temptablewidth}{1.0pt}}
\end{center}
\caption{\label{accu-3d-1} Accuracy test: the $L^1$ and $L^2$ errors and orders of accuracy for the three-dimensional advection of density perturbation obtained by linear CGKS-5 at $t = 2$.}
\begin{center}
\def\temptablewidth{0.8\textwidth}{\rule{\temptablewidth}{1.0pt}}
\begin{tabular*}{\temptablewidth}{@{\extracolsep{\fill}}c|c|c|c|c} 
MESH            & $L^1$ Error & Order   & $L^2$ Error & Order  \\
\hline
$10^3$ & $3.2325\times 10^{-2}$ &  & $1.2403\times 10^{-2}$  & \\
$20^3$ & $7.1472\times 10^{-4}$ & 5.4991 & $2.8195\times 10^{-4}$  & 5.4591 \\
$40^3$ & $2.3127\times 10^{-5}$ & 4.9497 & $9.0788\times 10^{-6}$  & 4.9568 \\
$80^3$ & $8.0092\times 10^{-7}$ & 4.8518 & $3.1656\times 10^{-7}$  & 4.8420 \\
\end{tabular*}
{\rule{\temptablewidth}{1.0pt}}
\end{center}
\caption{\label{accu-3d-2} Accuracy test: the $L^1$ and $L^2$ errors and orders of accuracy for the three-dimensional advection of density perturbation obtained by nonlinear CGKS-5 at $t = 2$.}
\end{table}

\subsection{Accuracy test}

In this subsection, the three-dimensional advection of density perturbation is used to validate the order of accuracy of the current compact scheme. 
The initial condition is given as follows
\begin{align*}
\rho_0&(x, y, z)=1+0.2\sin(\pi(x+y+z)),~p_0(x,y,z)=1,\\
&U_0(x,y,z)=1,~V_0(x,y,z)=1,~W_0(x,y,z)=1.
\end{align*}
The periodic boundary conditions are applied on all domain boundaries, and the
exact solution is
\begin{align*}
\rho(x,y&,z,t)=1+0.2\sin(\pi(x+y+z-3t)),~p(x,y,z,t)=1,\\
&U(x,y,z,t)=1,~V(x,y,z,t)=1,~W(x,y,z,t)=1.
\end{align*}
The computational domain is $[0,2]\times[0,2]\times[0,2]$. 
In this case,  $\tau=0$ and the gas-distribution function reduces to
\begin{align*}
f(\boldsymbol{x}_{G},t,\boldsymbol{u},\xi)=g(1+At).
\end{align*}
The $L^1$ errors, $L^2$ errors and orders of accuracy obtained by the current scheme with linear and nonlinear reconstructions at $t=2$ are presented in Tab.\ref{accu-3d-1} and Tab.\ref{accu-3d-2}, respectively. 
It can be observed that both the linear and nonlinear CGKS-5 achieve the expected fifth-order of accuracy.
Except for the results obtained on the coarsest mesh, the nonlinear scheme demonstrates consistent errors and order of accuracy when compared to the linear scheme.

\begin{figure}[!h]
\centering
\includegraphics[width=0.485\textwidth]{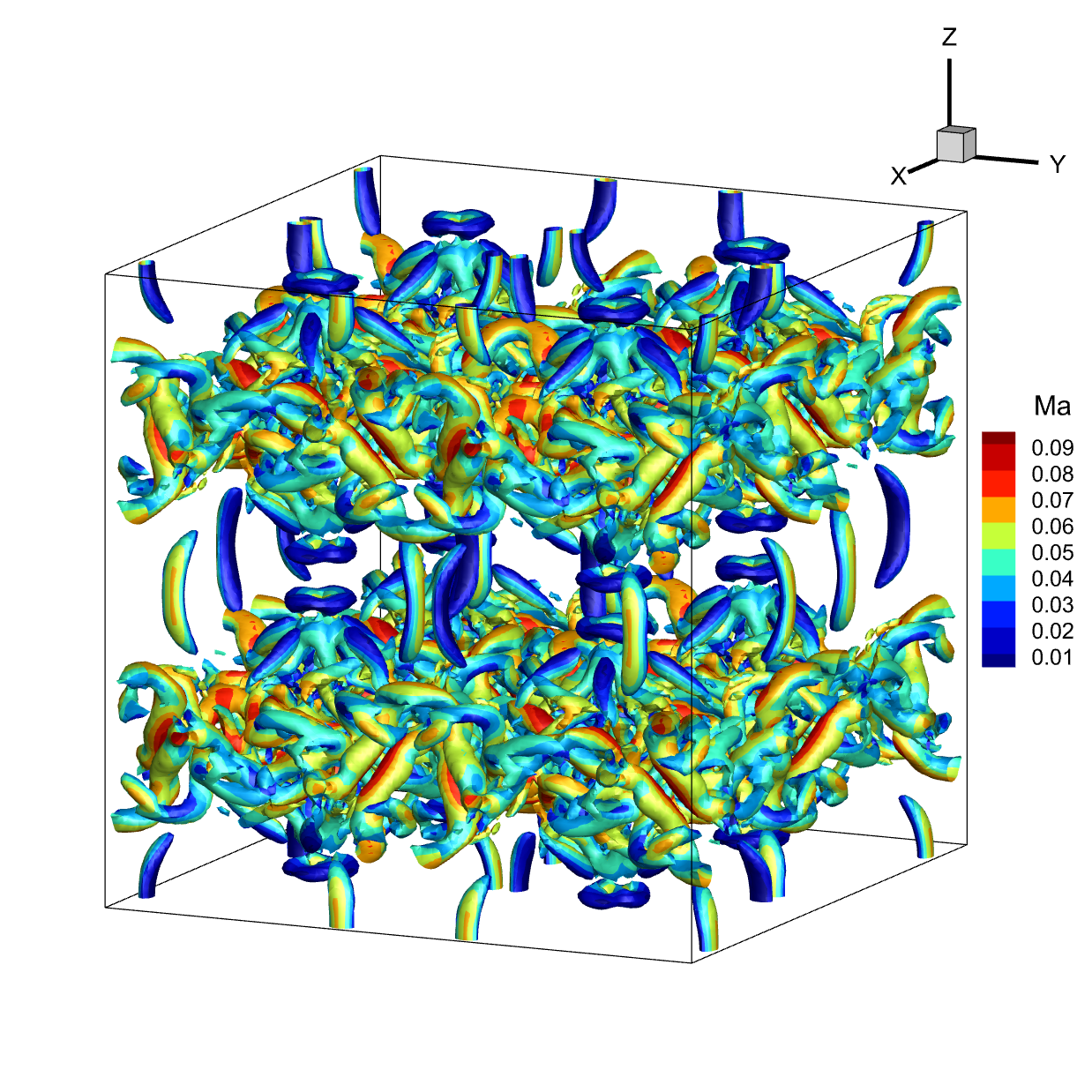}
\includegraphics[width=0.485\textwidth]{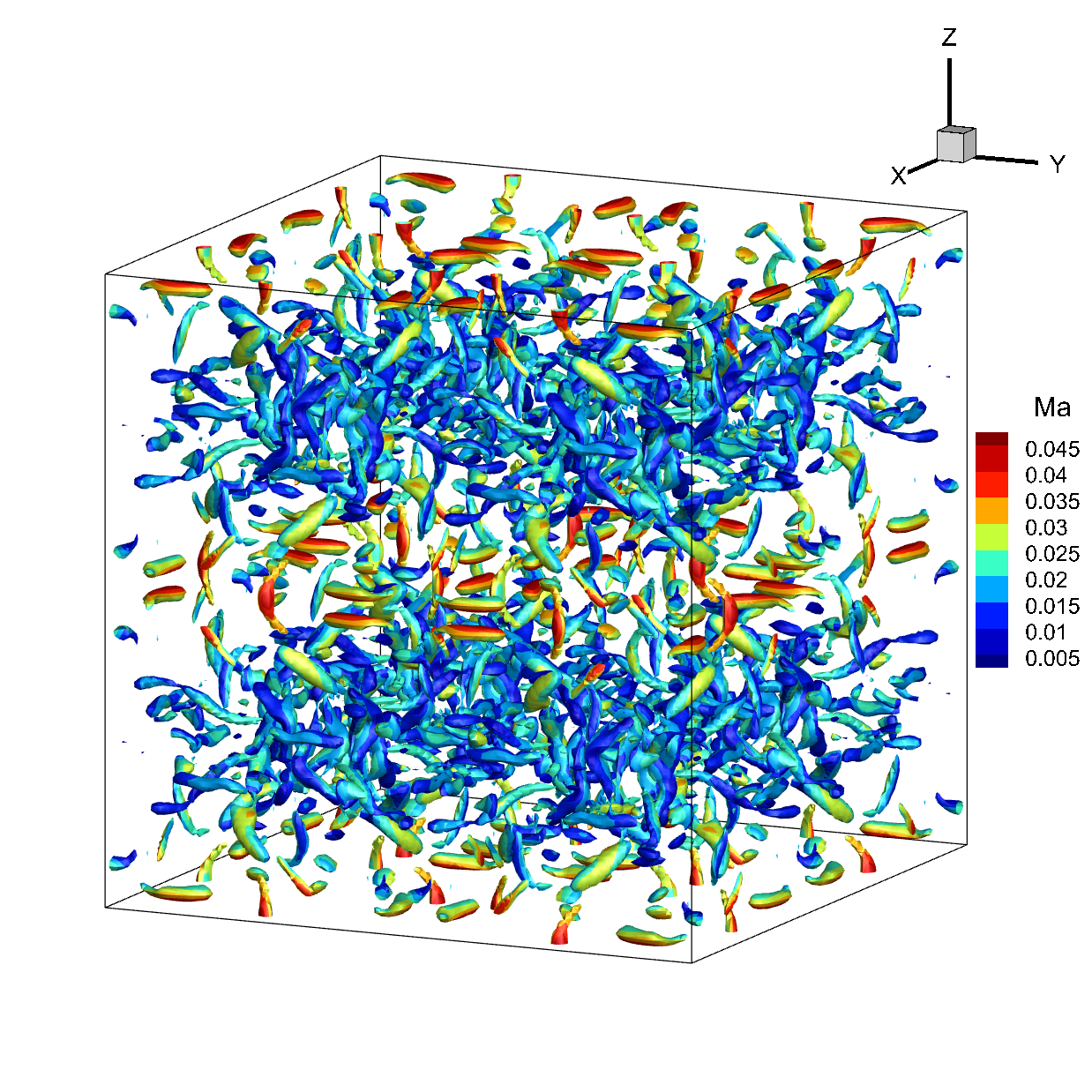}
\caption{\label{TGV-sub-1} Taylor-Green vortex flow: Iso-surfaces of the Q-criterion ($Q=1.5$) at $t=10$ (left) and $t=20$ (right) for $Ma_{\infty}=0.1$ on a $128^3$ mesh.}
\includegraphics[width=0.485\textwidth]{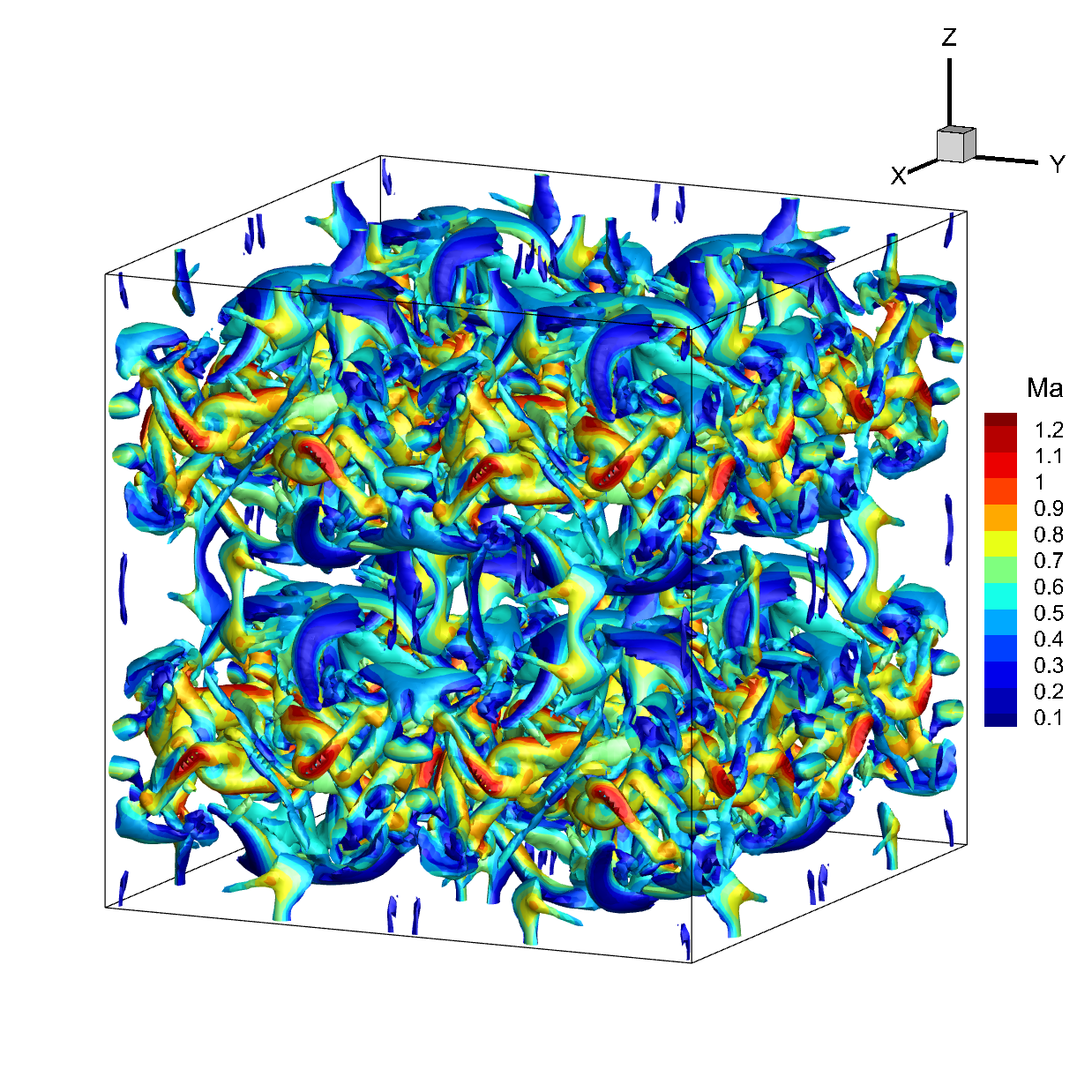}
\includegraphics[width=0.485\textwidth]{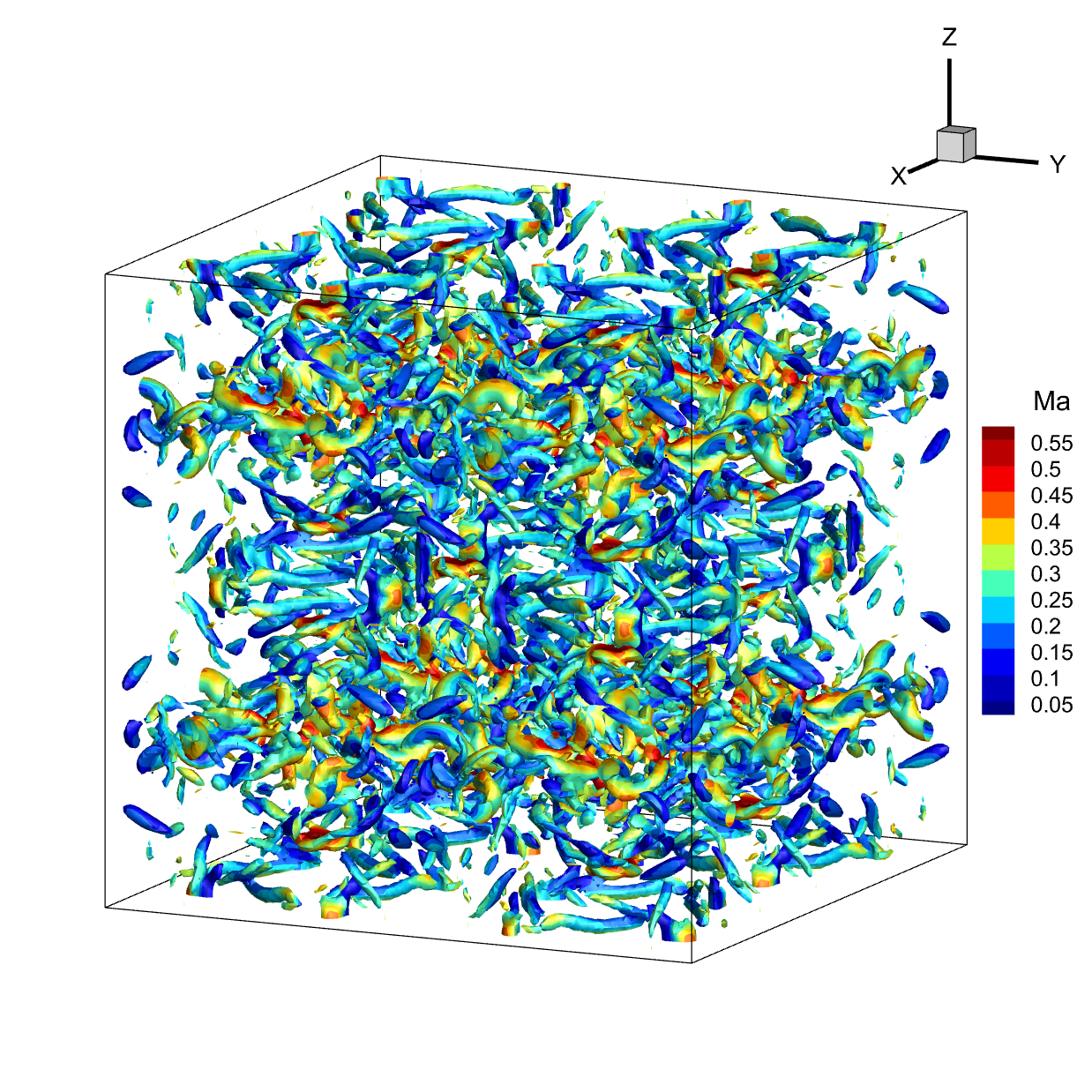}
\caption{\label{TGV-sup-1} Taylor-Green vortex flow: Iso-surfaces of the Q-criterion ($Q=1.5$) at $t=10$ (left) and $t=20$ (right) for $Ma_{\infty}=1.25$ on a $128^3$ mesh.}
\end{figure}

\begin{figure}[!h]
\centering
\includegraphics[width=0.485\textwidth]{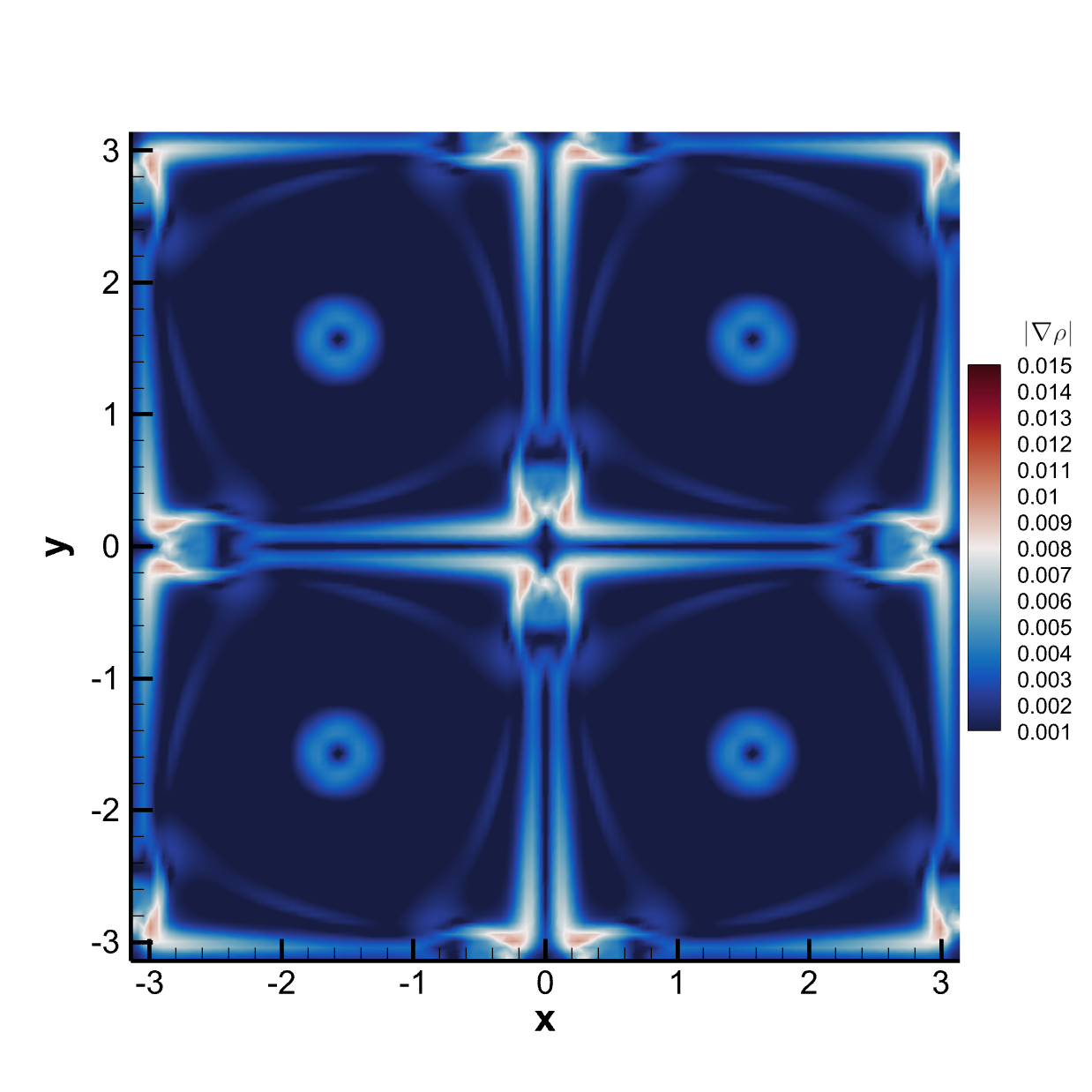}
\includegraphics[width=0.485\textwidth]{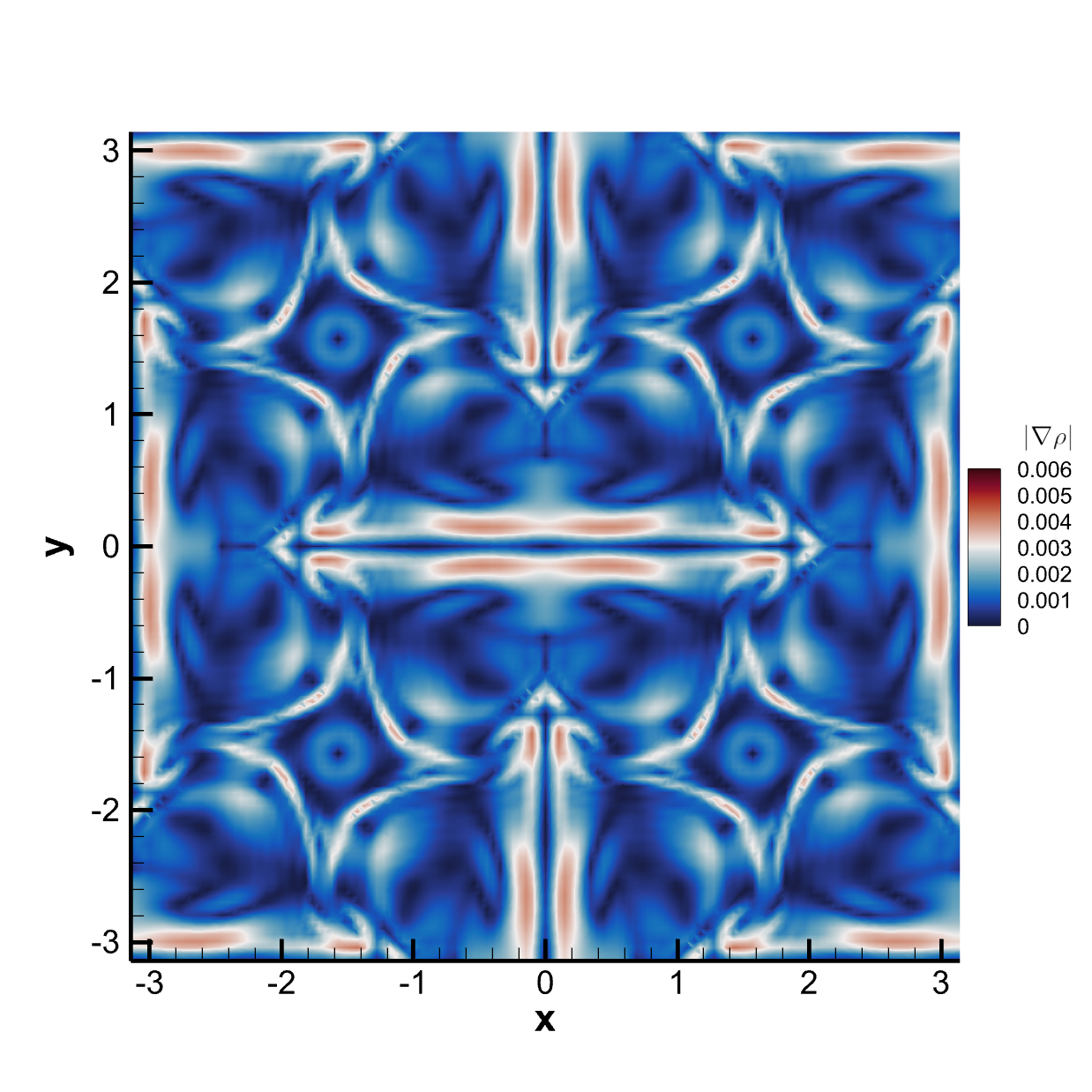}
\caption{\label{TGV-sub-2} Taylor-Green vortex flow: the magnitude of the density gradient $\nabla\rho$ for $Ma_{\infty}=0.1$ at $t=10$ (left) and $t=20$ (right) on the centerline (x-y) plane at $z=0$.}
\includegraphics[width=0.485\textwidth]{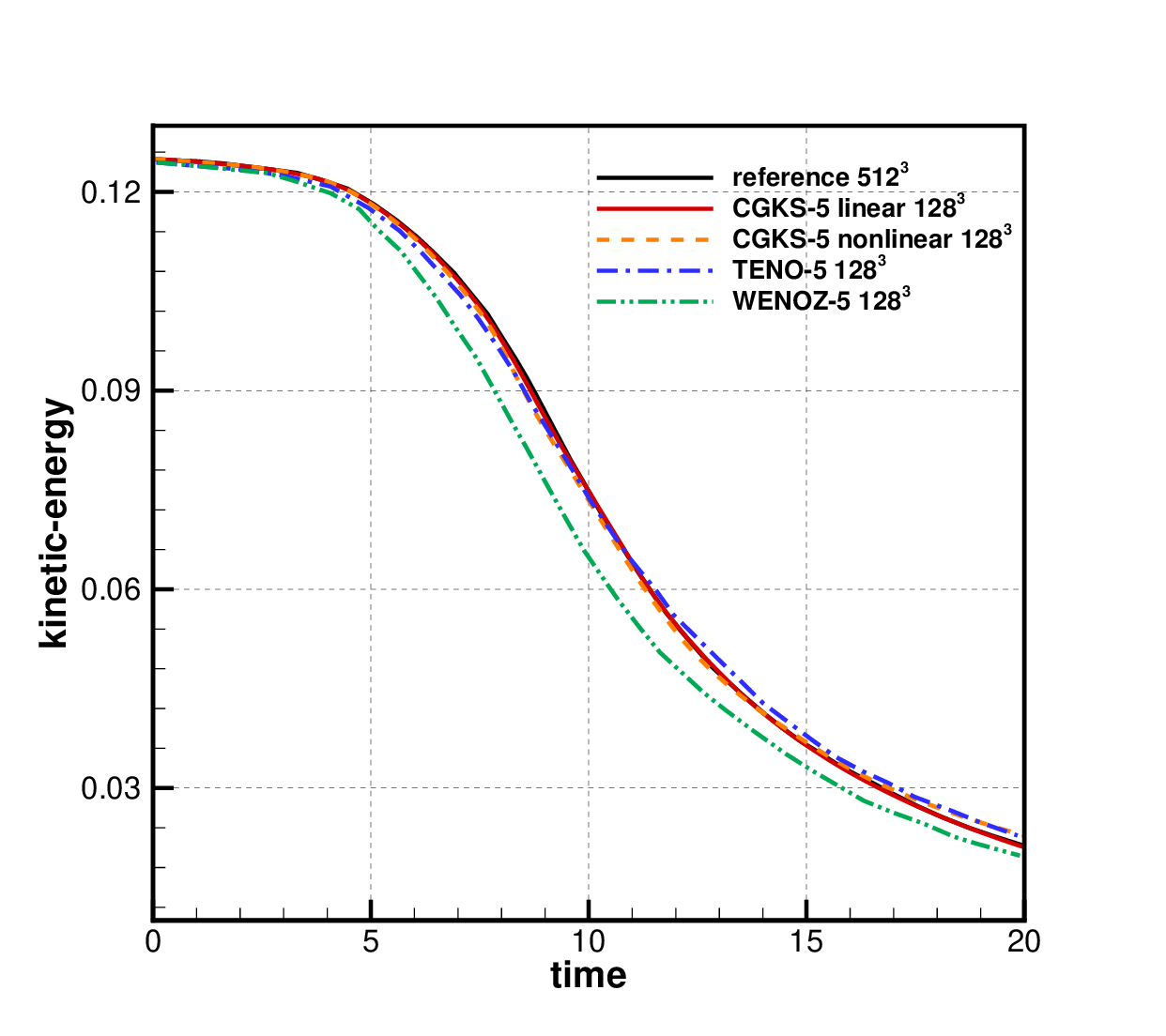}
\includegraphics[width=0.485\textwidth]{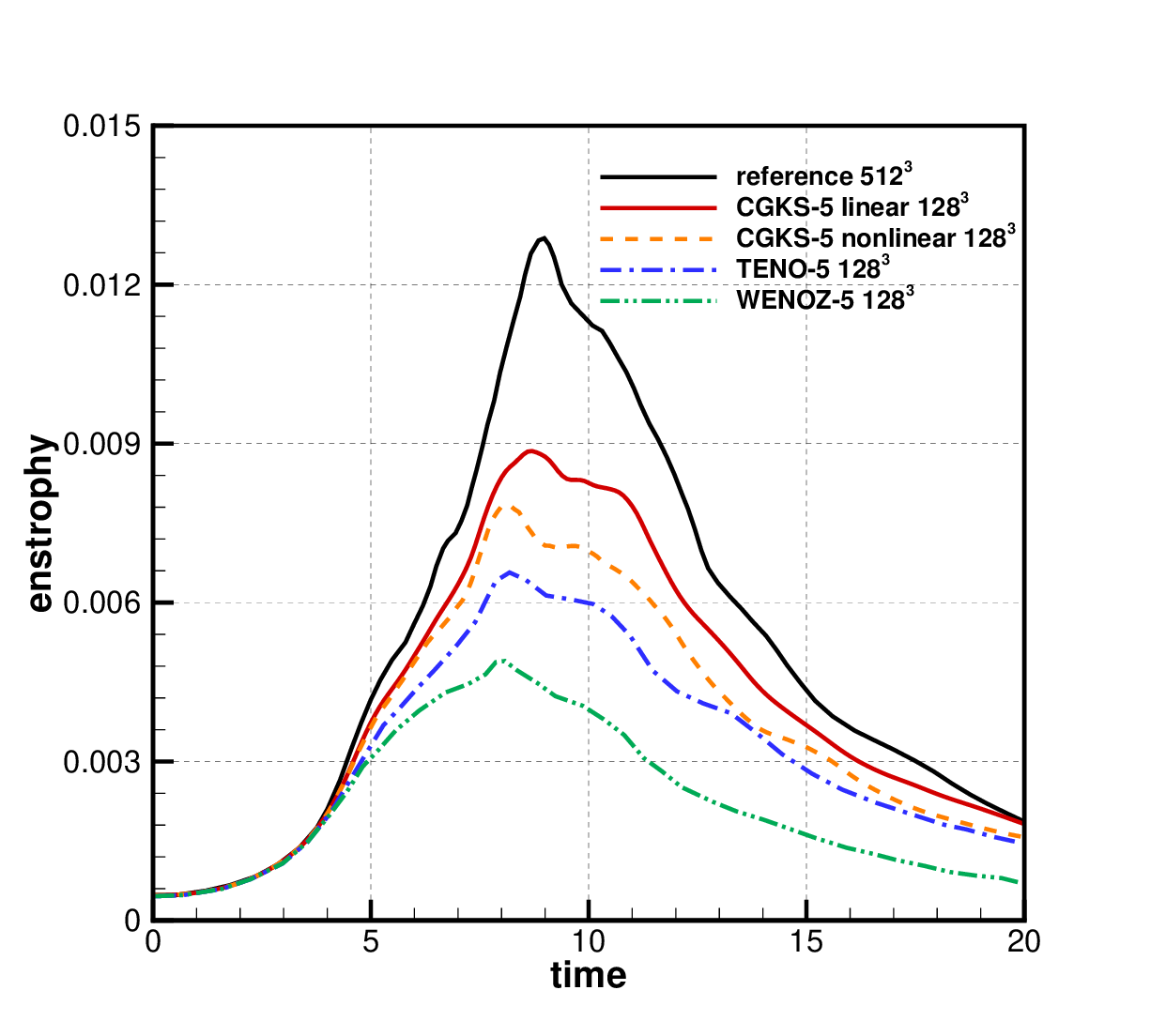}
\caption{\label{TGV-sub-3} Taylor-Green vortex flow: the comparison of the time histories of kinetic energy (left) and enstrophy (right) using different high-order schemes for $Ma_{\infty}=0.1$.}
\end{figure}

\begin{figure}[!h]
\centering
\includegraphics[width=0.485\textwidth]{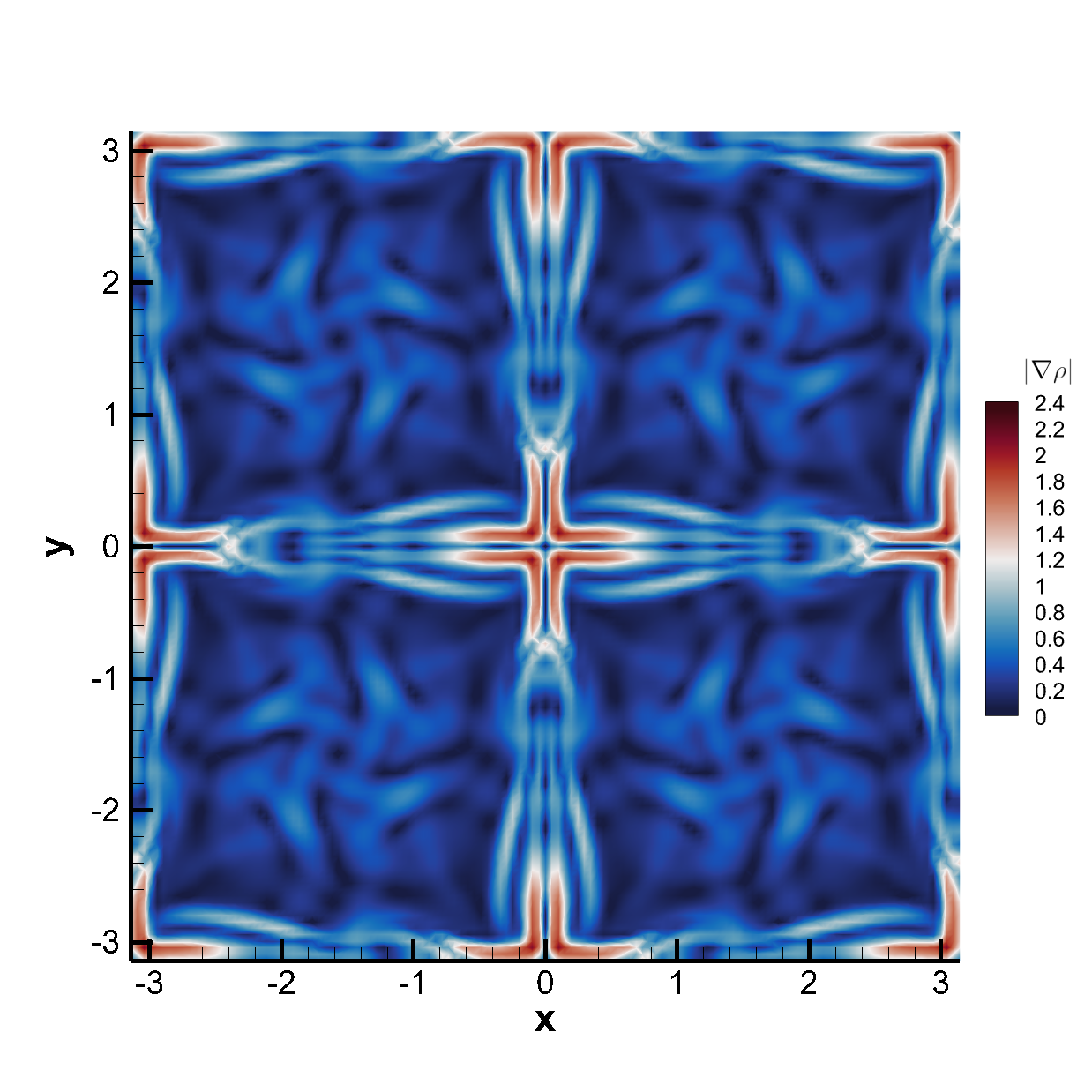}
\includegraphics[width=0.485\textwidth]{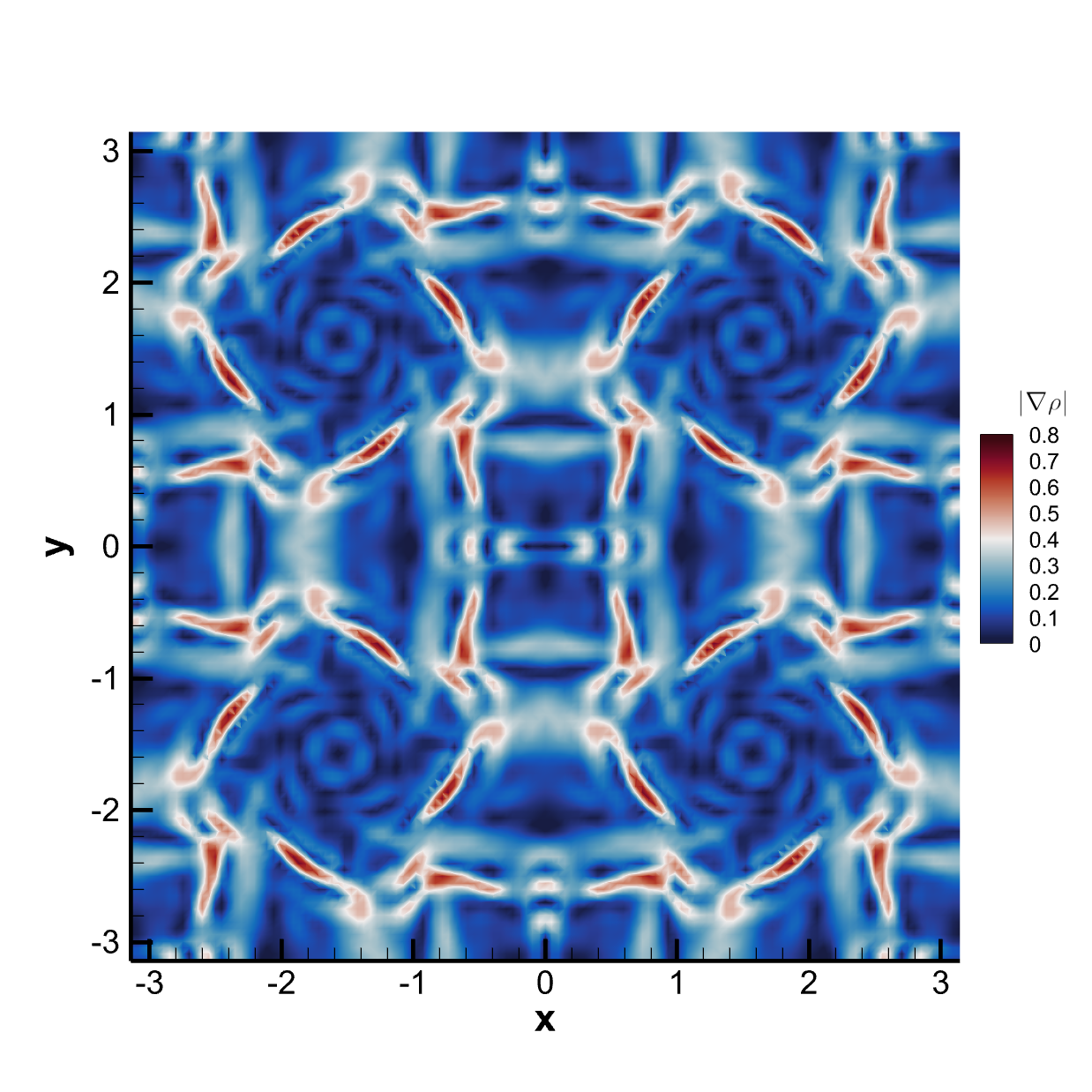}
\caption{\label{TGV-sup-2} Taylor-Green vortex flow: the magnitude of the density gradient $\nabla\rho$ for $Ma_{\infty}=1.25$ at $t=10$ (left) and $t=20$ (right) on the centerline (x-y) plane at $z=0$.}
\includegraphics[width=0.485\textwidth]{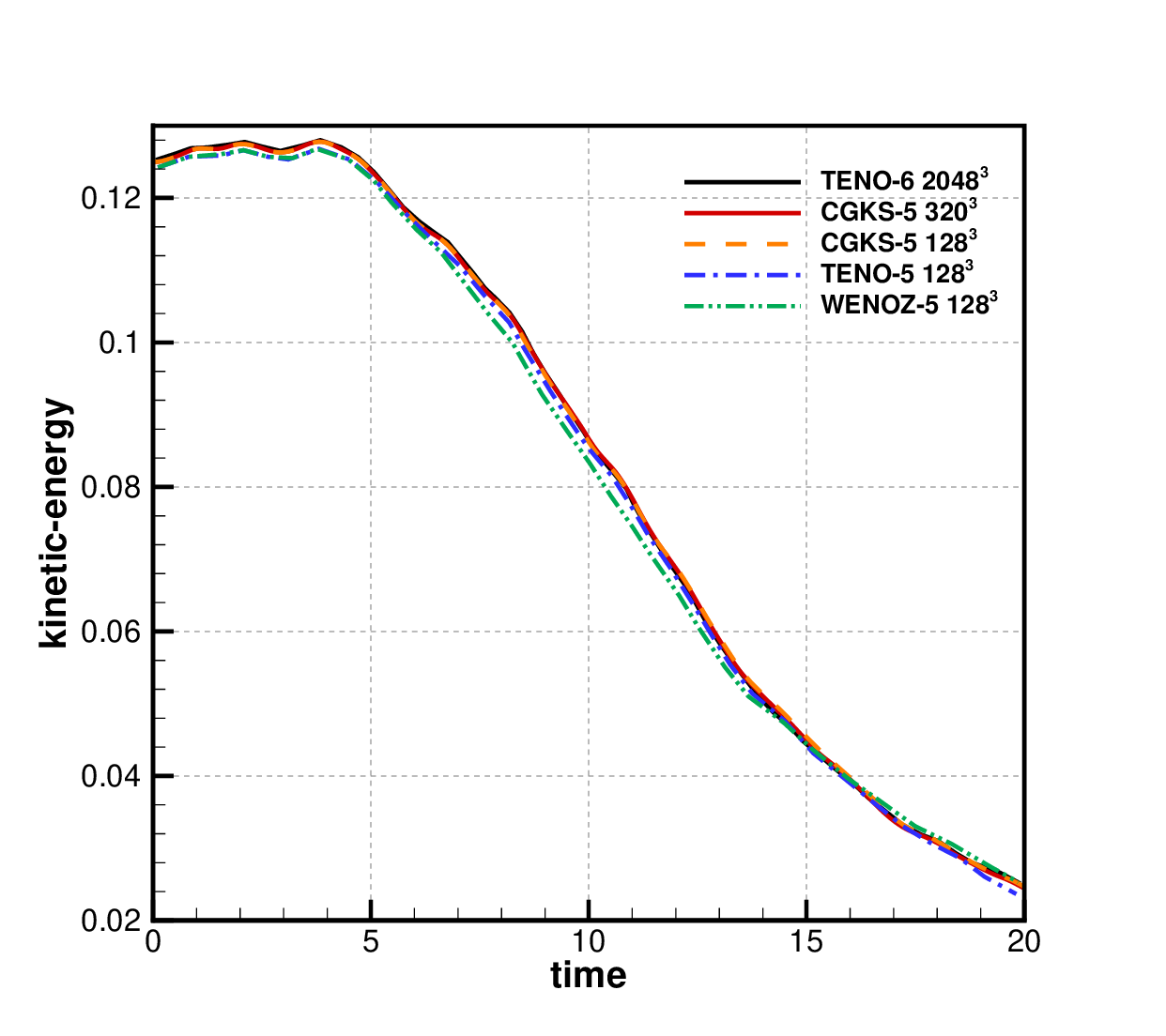}
\includegraphics[width=0.485\textwidth]{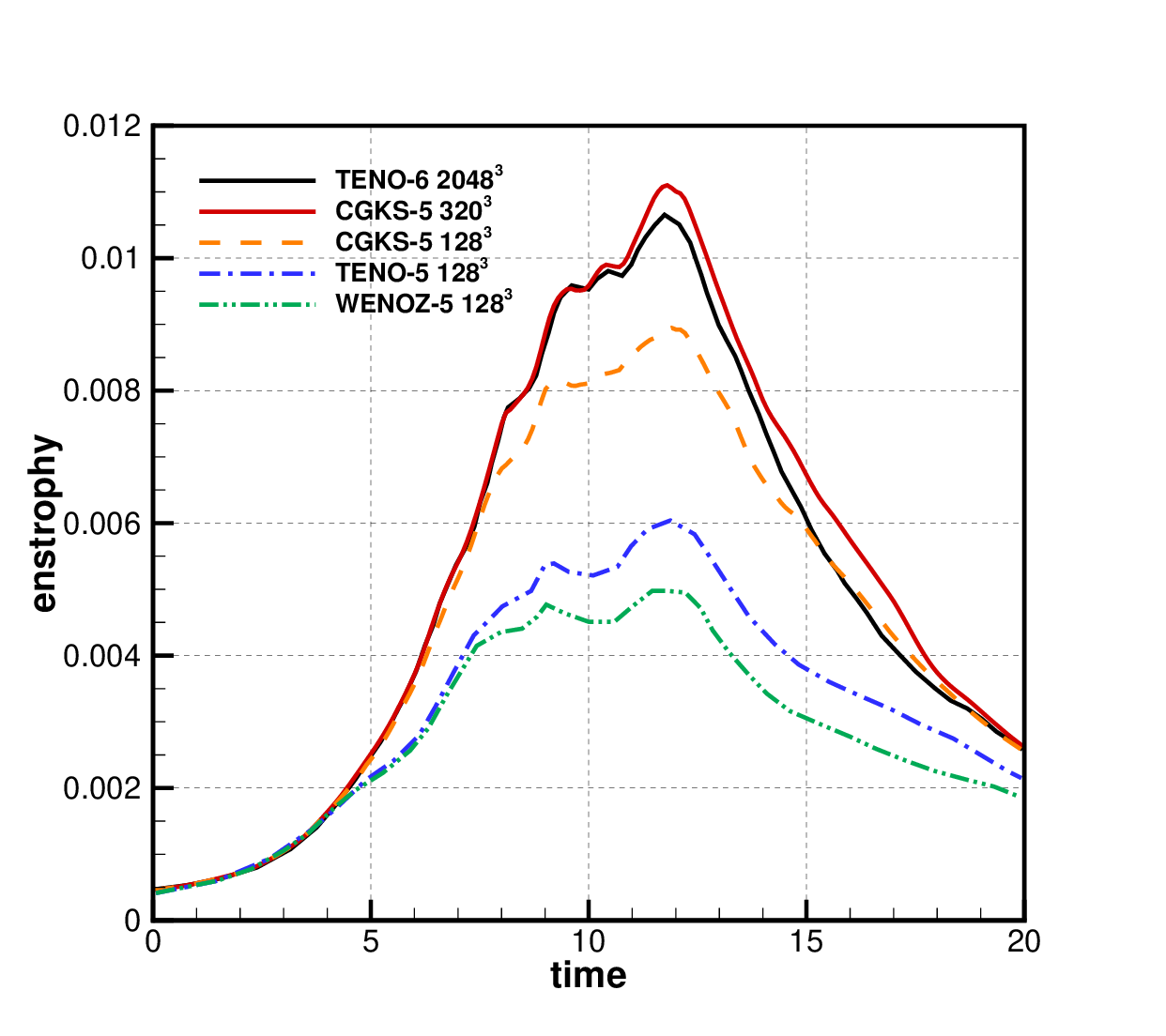}
\caption{\label{TGV-sup-3} Taylor-Green vortex flow: the comparison of the time histories of kinetic energy (left) and enstrophy (right) using different high-order schemes for $Ma_{\infty}=1.25$.}
\end{figure}

\begin{figure}[!h]
\centering
\includegraphics[width=0.485\textwidth]{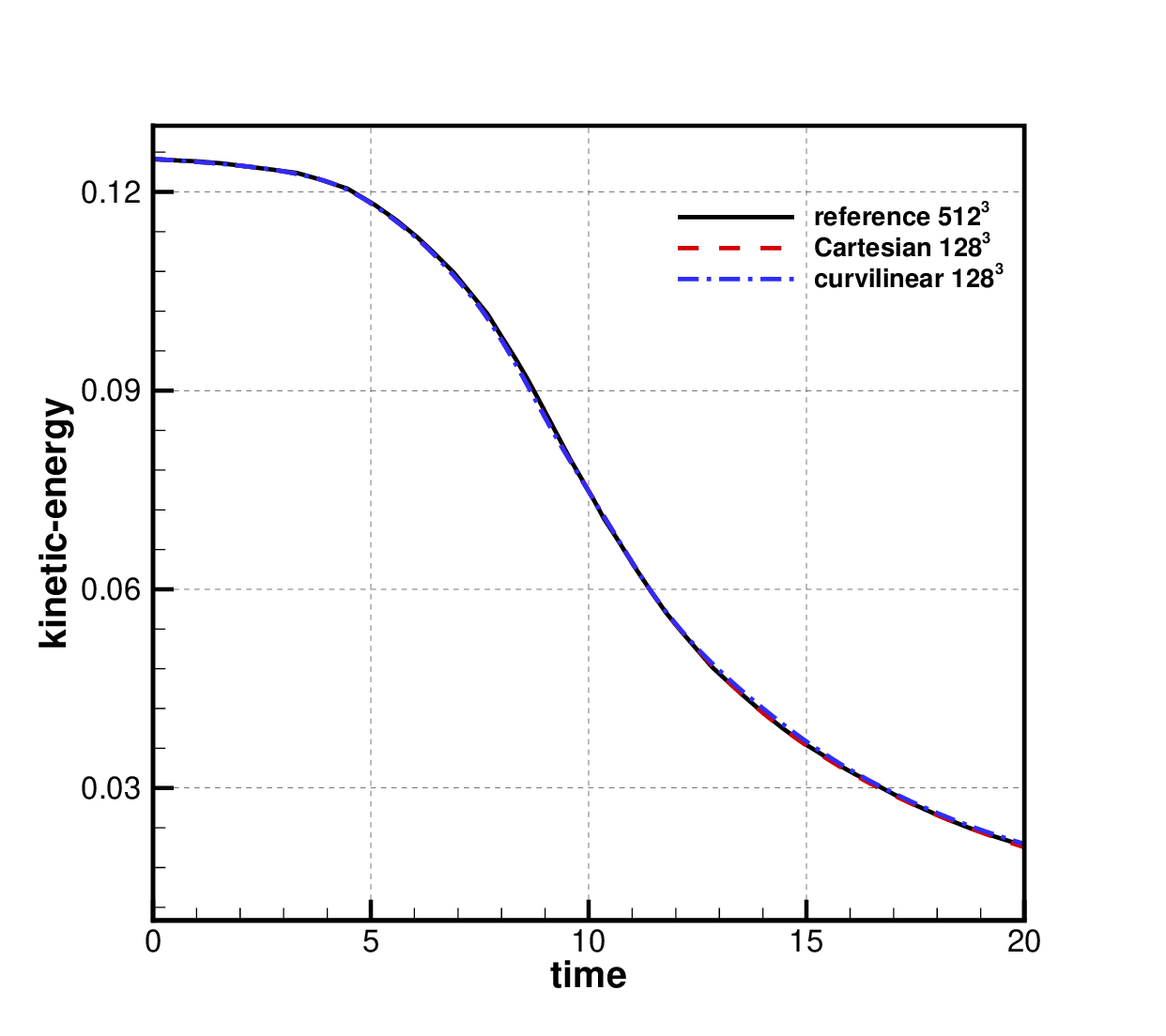}
\includegraphics[width=0.485\textwidth]{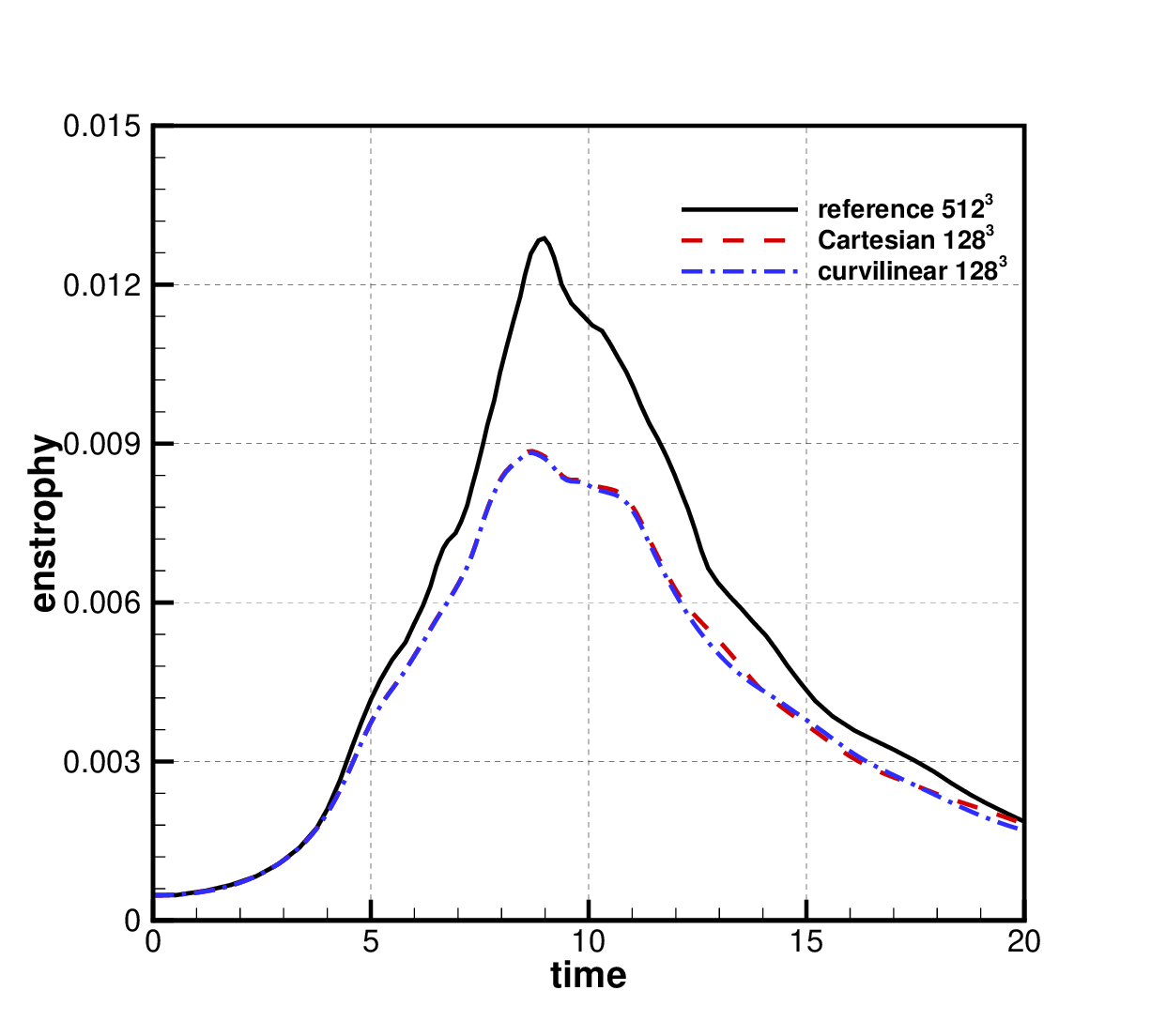}
\caption{\label{TGV-curv} Taylor-Green vortex flow:  the comparison of the time histories of kinetic energy (left) and enstrophy (right) using Cartesian and curvilinear meshes for $Ma_{\infty}=0.1$.}
\end{figure}

\subsection{Taylor-Green vortex flow}

The Taylor-Green vortex flow is a well-known problem in computational fluid dynamics used to study vortex dynamics, turbulent decay and energy cascade mechanisms \cite{TGV-1,TGV-2,TGV-3}. It is particularly useful for verifying the accuracy and stability of the current fifth-order compact scheme for viscous flow. 
The initial velocity field is given by
\begin{align*}
U &= U_0\sin(\frac{x}{L})\cos(\frac{y}{L})\cos(\frac{z}{L}), \\
V &=-U_0\cos(\frac{x}{L})\sin(\frac{y}{L})\cos(\frac{z}{L}), \\
W &= 0, 
\end{align*}
and the initial pressure field is set as
\begin{align*}
\rho &= \rho_0\big(1+\frac{\gamma Ma_{\infty}^2}{16}\big(\cos(\frac{2x}{L})+\cos(\frac{2y}{L})\big)\big(\cos(\frac{2z}{L})+2\big)\big),\\
p &= p_0\big(1+\frac{\gamma Ma_{\infty}^2}{16}\big(\cos(\frac{2x}{L})+\cos(\frac{2y}{L})\big)\big(\cos(\frac{2z}{L})+2\big)\big).
\end{align*}
In the computation, $L=1$, $U_0=1$, $\rho_0=1$, and the Prandtl number is set as 0.7.
The Reynolds number is $Re=1600$ defined by $Re=\rho_0U_0L/\mu_0$ , where 
$\mu_0$ is the dynamic viscosity coefficient. 
The Mach number is determined by $Ma_{\infty}=U_0/\sqrt{\gamma p_0/\rho_0}$.
This case is conducted at two Mach numbers, $Ma_{\infty}=0.1$,  which falls within the incompressible limit, and $Ma_{\infty}=1.25$, where shock waves are present.
The Taylor-Green vortex flow is defined in a periodic domain defined as $[-\pi L, \pi L]^3$.
In the subsonic case, both linear and nonlinear reconstructions are used, while in the supersonic case, only nonlinear reconstruction is employed.
A uniform $128^3$ mesh is used for the two configurations.

Figure \ref{TGV-sub-1} and Figure \ref{TGV-sup-1} present the iso-surfaces of the Q-criterion ($Q=1.5$), colored according to the Mach number, at $t=10$ and $t=20$.  These figures indicate the decay of vortices over time and demonstrate the energy transfer from larger scales to smaller scales.
Figure \ref{TGV-sub-2} and Figure \ref{TGV-sup-2} present the magnitude of the density gradient for $Ma_{\infty}=0.1$ and $Ma_{\infty}=1.25$ at $t=10$ and $t=20$, respectively, on the centerline (x-y) plane at $z=0$.  
These figures demonstrate that the flow structures can be accurately resolved with high resolution using the present compact scheme for both subsonic and supersonic cases.

To evaluate the performance of CGKS-5 more accurately, several averaged quantities are calculated as it evolves over time.
The volume-averaged kinetic energy is given by
\begin{equation*}
E_k=\frac{1}{\rho_0 |\Omega|}\int_{\Omega} \frac{1}{2}\rho \boldsymbol{U}\cdot\boldsymbol{U}\,\mathrm{d}V,
\end{equation*}
where $\Omega$ is the computational domain and $E_k$ is calculated by numerical quadrature. The integrated enstrophy is defined as
\begin{equation*}
\epsilon^S=\frac{1}{\rho_0 |\Omega|}\int_{\Omega}\mu(|\nabla\times\boldsymbol{U}|^2)\,\mathrm{d}V,
\end{equation*}
where the dynamic viscosity $\mu$ is computed by the Sutherland's law, $\epsilon^S$ is calculated by numerical quadrature, and the velocity derivative values at quadrature points are obtained by the compact reconstruction.

The comparison of the time histories of kinetic energy $E_k$ and enstrophy $\epsilon^S$ using different high-order schemes is presented in Figure \ref{TGV-sub-3} and Figure \ref{TGV-sup-3}.
The results obtained using the fifth-order WENO-Z and fifth-order TENO schemes, as shown in Figure \ref{TGV-sub-3} and Figure \ref{TGV-sup-3}, are computed within the framework of the finite volume method on a $128^3$ mesh \cite{TGV-ref-2}.
The reference line in Figure \ref{TGV-sub-3} represents a spectral solution computed on a $512^3$ mesh \cite{TGV-ref-1}, while the result obtained by the sixth-order TENO scheme in Figure \ref{TGV-sup-3} is calculated on a $2048^3$ mesh \cite{TGV-ref-3}.
For the subsonic case, Figure \ref{TGV-sub-3} shows that the results calculated by CGKS-5, using both linear and nonlinear reconstructions, closely align with the reference solution for kinetic energy. At the same time, the result obtained by the nonlinear CGKS-5 demonstrates higher resolution in enstrophy dissipation compared to the nonlinear WENO-Z and TENO schemes.
For the supersonic case, Figure \ref{TGV-sup-3} shows that CGKS-5 achieves significantly higher resolution in enstrophy dissipation. Additionally, it can be observed that the current compact scheme exhibits excellent mesh convergence as the resolution increases from the $128^3$ mesh to the $320^3$ mesh.

To verify the stability and resolution of CGKS-5 on curvilinear mesh, 
the subsonic case for $Ma_{\infty}=0.1$ and $Re=1600$ with linear reconstruction is tested.
The following slight distortion is applied to the uniform Cartesian mesh
\begin{equation*}
\begin{cases}
\hat{x}=x+C\cdot \sin x\sin y,\\
\hat{y}=y+C\cdot \sin x\sin y,\\
\hat{z}=z,
\end{cases}
\end{equation*}
where $(x,y,z)\in [-\pi L,\pi L]^3$, $\Delta x=\Delta y=\Delta z=2\pi L/N$ and $C=0.1$.
Figure \ref{TGV-curv} presents the comparison of the time histories of kinetic energy $E_k$ and enstrophy $\epsilon^S$ using Cartesian and curvilinear meshes.
It can be observed that the results on the curvilinear mesh are consistent with those on the orthogonal uniform mesh, maintaining high resolution and linear stability.

\begin{figure}[!h]
\centering
\includegraphics[width=0.47\textwidth]{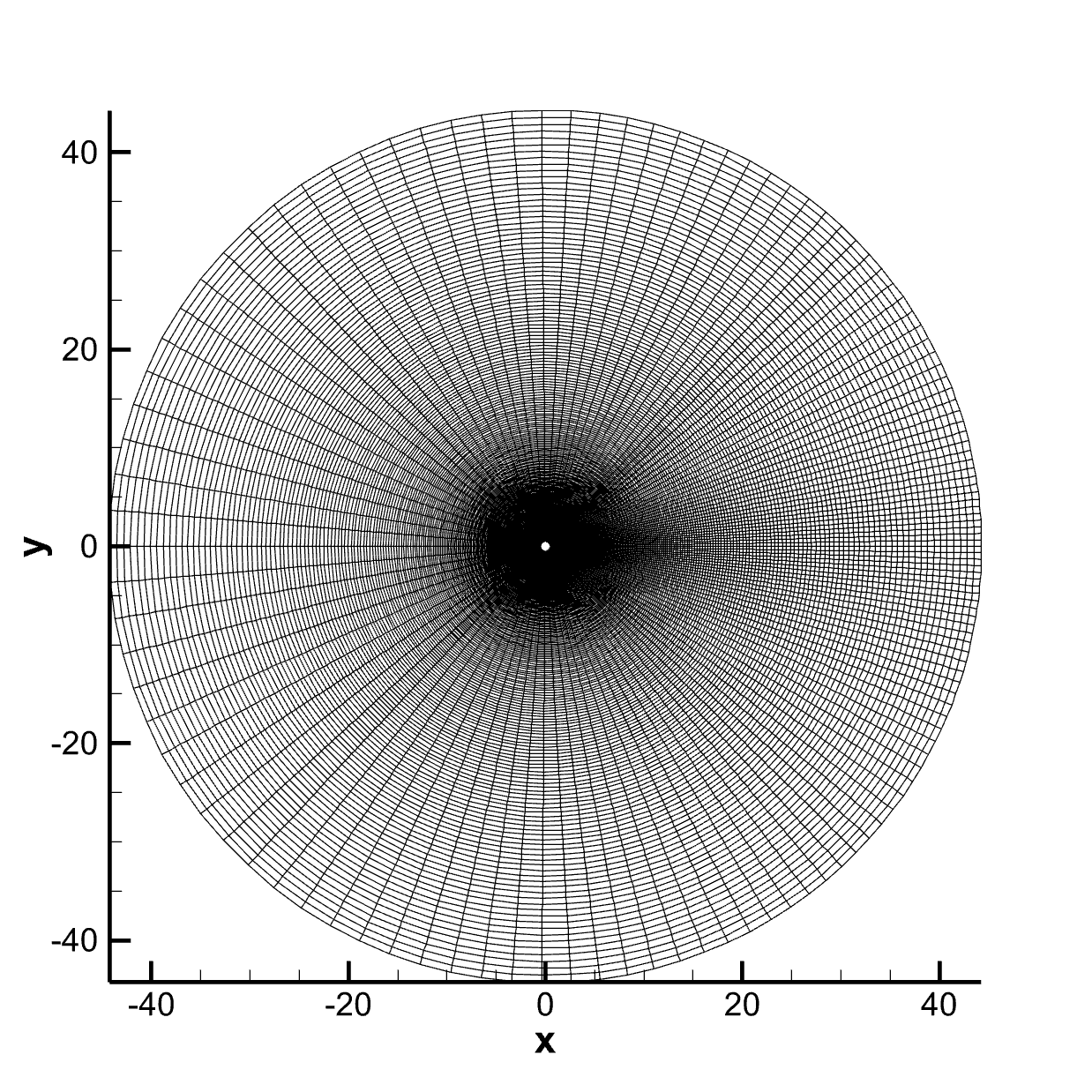}
\includegraphics[width=0.47\textwidth]{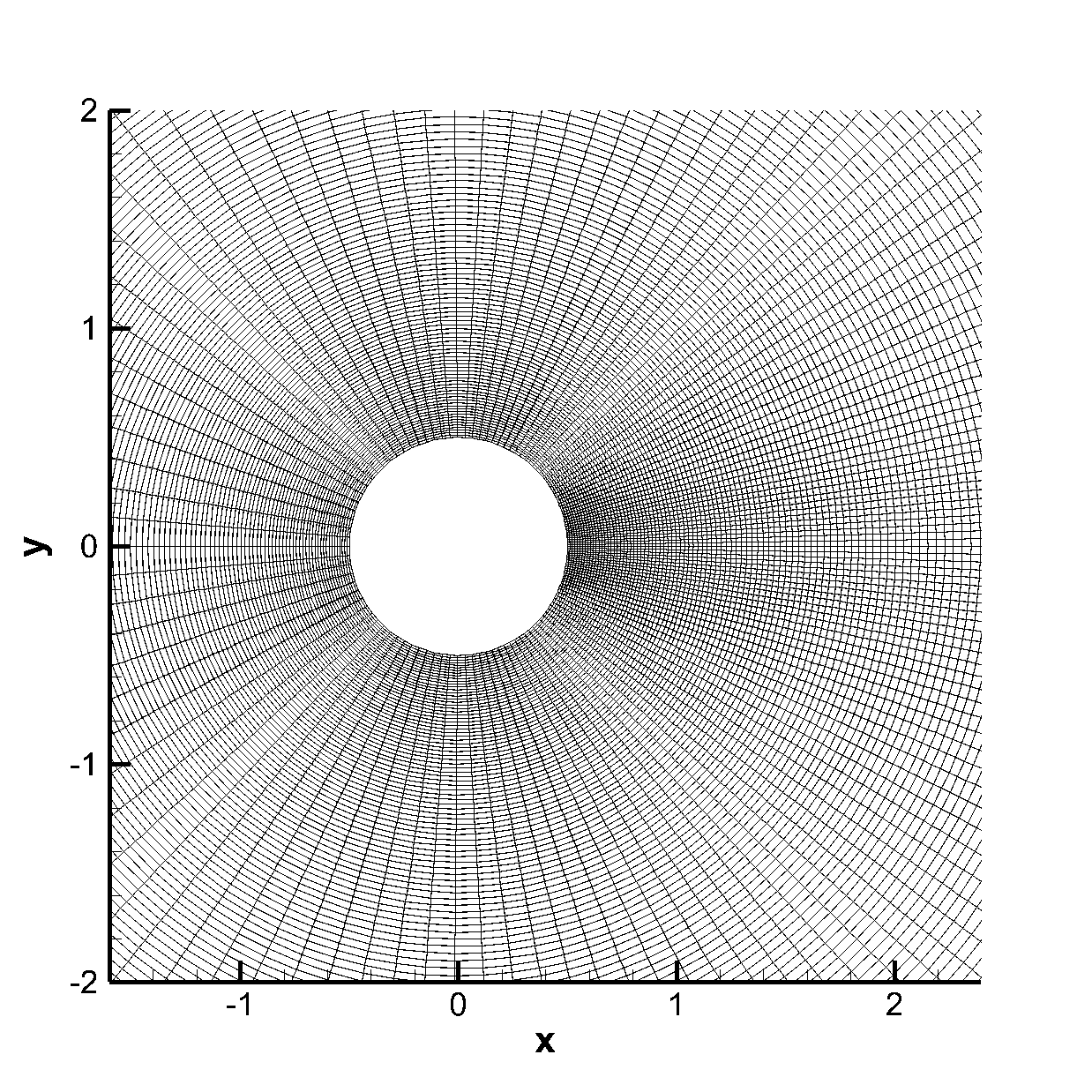}
\caption{\label{cylinder-mesh} Turbulent flow past a cylinder: the visualization of the computational mesh on the x-y plane with a full view (left) and a detailed view (right).}
\end{figure}

\begin{figure}[!h]
\centering
\includegraphics[width=0.49\textwidth]{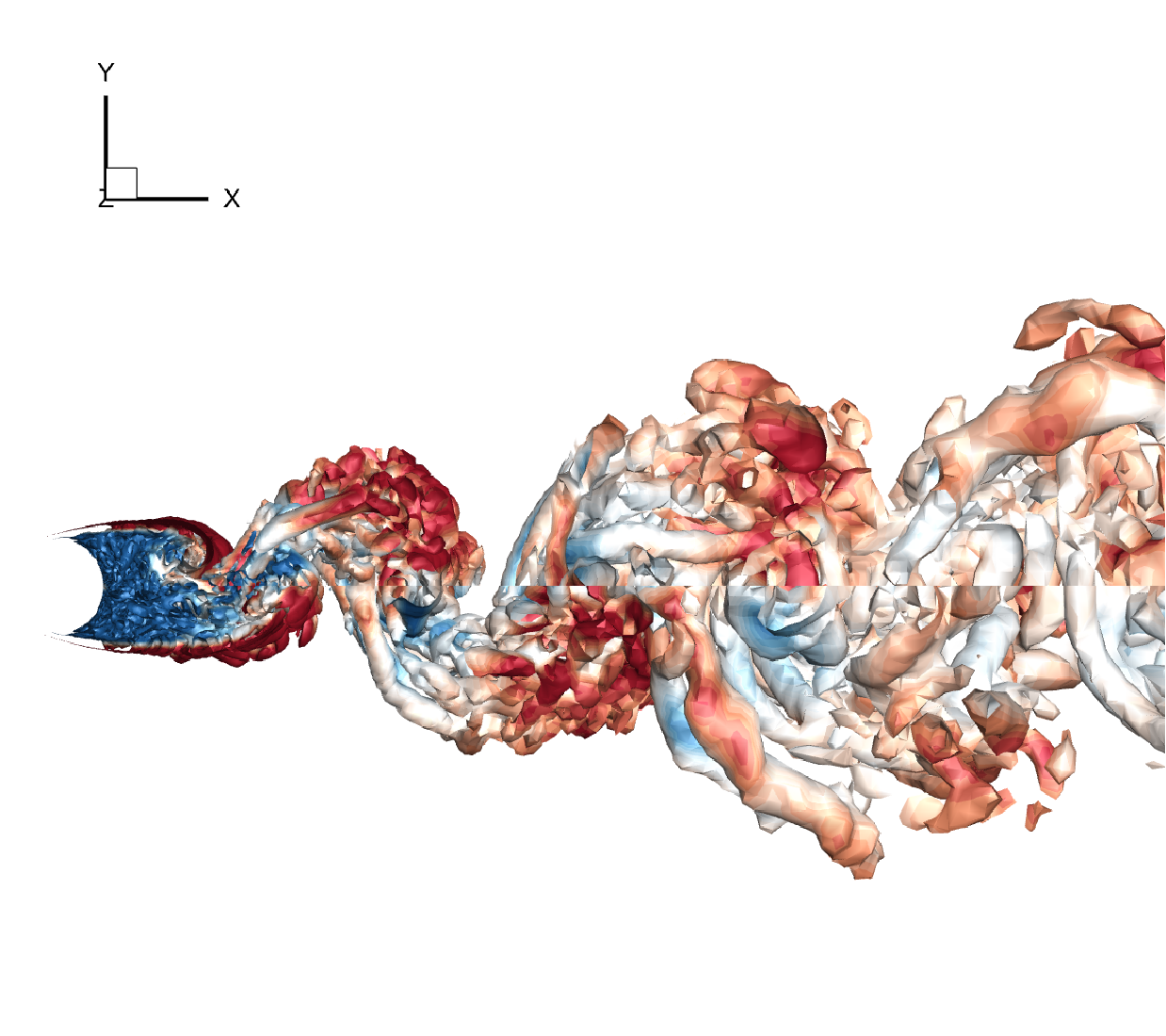}
\includegraphics[width=0.485\textwidth]{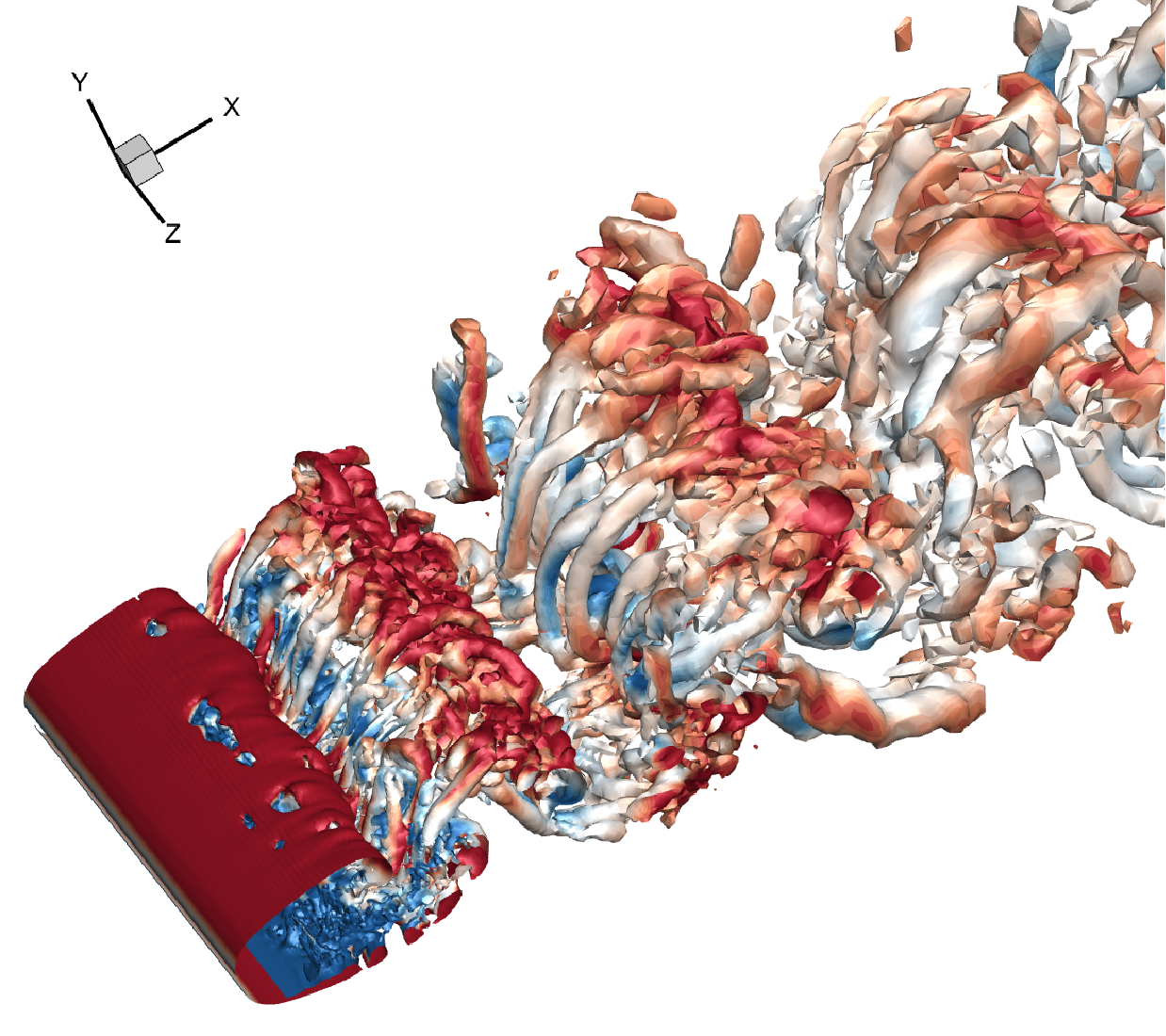}
\caption{\label{cylinder-vortex} Turbulent flow past a cylinder: the iso-surface of the Q-criterion ($Q=0.01$) with a front view (left) and an isometric view (right).}
\end{figure}

\begin{figure}[!h]
\centering
\includegraphics[width=0.49\textwidth]{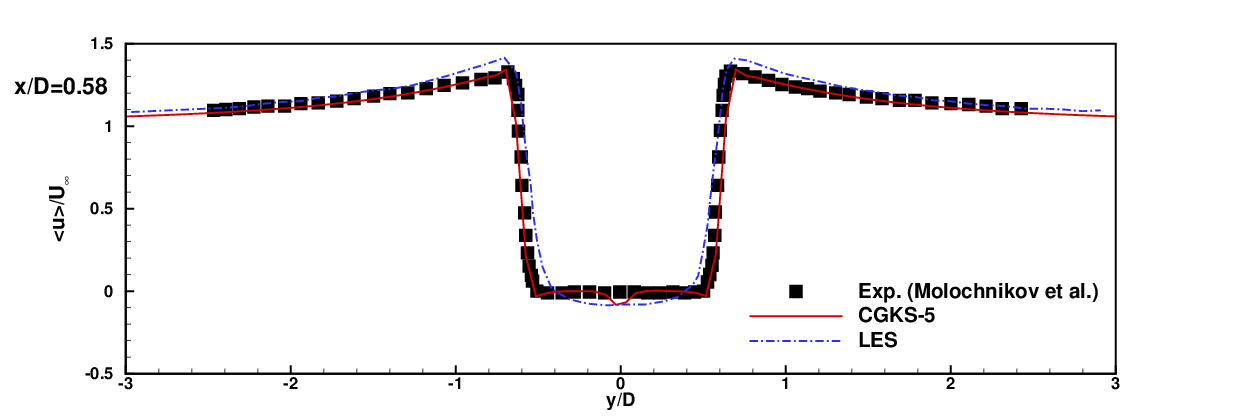}
\includegraphics[width=0.49\textwidth]{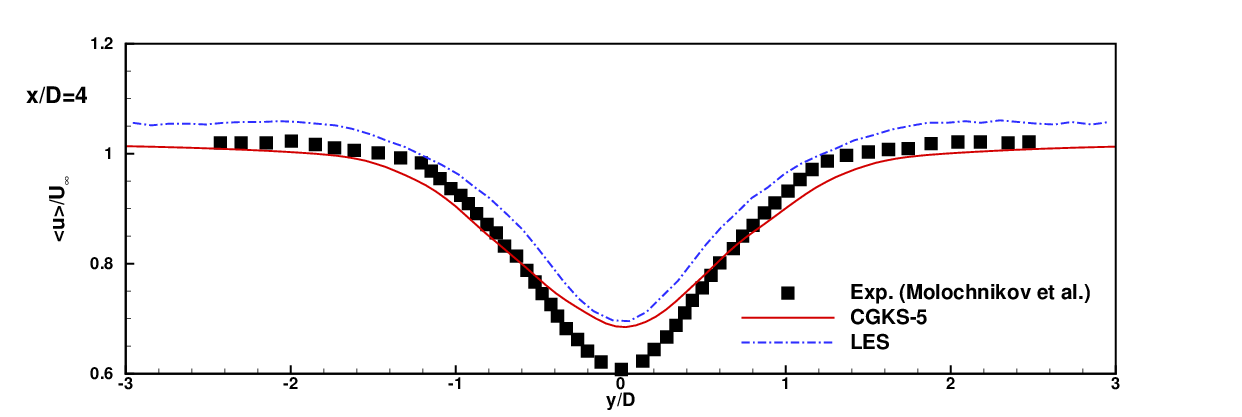}
\includegraphics[width=0.49\textwidth]{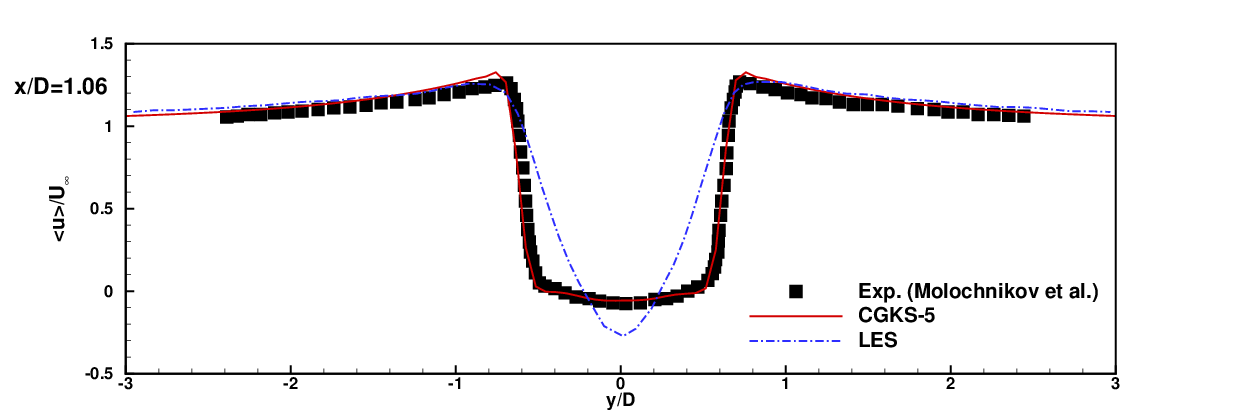}
\includegraphics[width=0.49\textwidth]{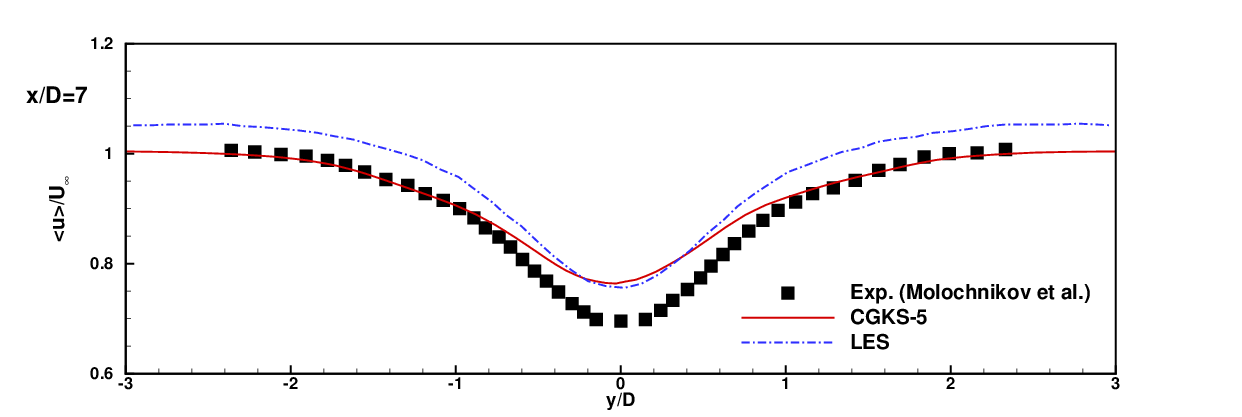}
\includegraphics[width=0.49\textwidth]{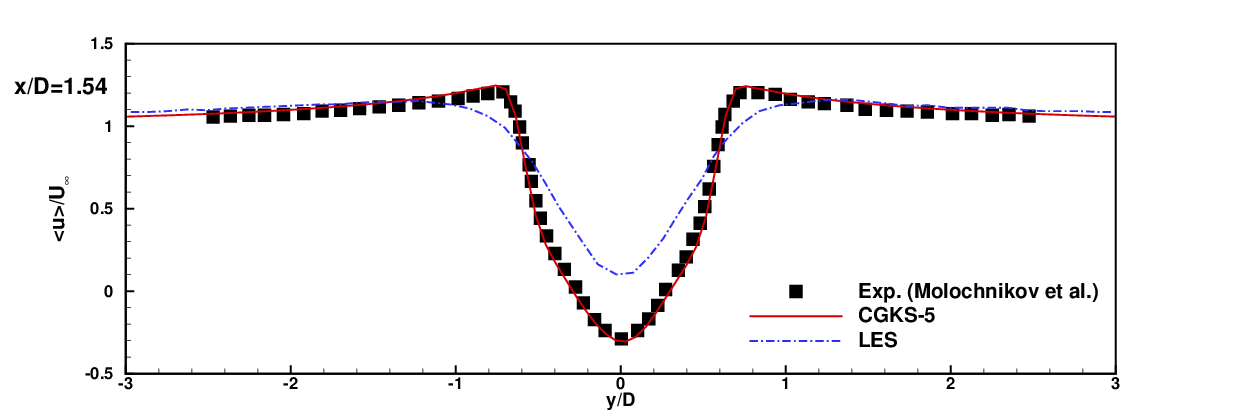}
\includegraphics[width=0.49\textwidth]{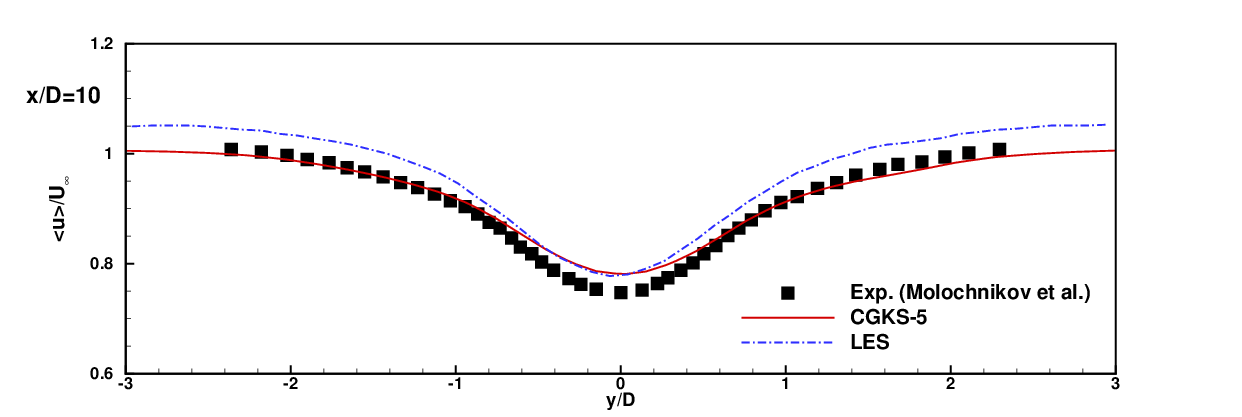}
\caption{\label{cylinder-line} Turbulent flow past a cylinder: the streamwise velocity profiles at different locations in the cylinder wake obatained by the linear fifth-order CGKS.}
\end{figure}

\subsection{Turbulent flow past a cylinder}

The numerical performance of the current compact scheme on curvilinear meshes is validated using the turbulent flow past a circular cylinder.
The cylinder, with a diameter of $D=1$, is positioned at the origin, and its wall defines the inner boundary of the physical domain.
The outer boundary of the domain is concentric with the cylinder, with a diameter of approximately $40D$.
This case is set in a three-dimensional domain, where the z-direction is periodic with a length of $L_z = 2\pi D$.

The initial condition is given as a free stream condition
\begin{align*}
(\rho,U,V,W,p)_{\infty} = (1, Ma_\infty,0,0, 1/\gamma),
\end{align*}
where $Ma_\infty=0.25$ is the Mach number of the free stream, and the Reynolds number $Re=3900$.
The non-slip adiabatic boundary condition is imposed on the wall of the cylinder, and the far-field boundary condition is applied to the outer boundary.
The computational mesh is an “O-” type structured mesh with $256\times 128\times 128$ cells, where the size of the first layer cells on the cylinder is $h_{min}=1.1\times 10^{-2}$ and the mesh spacing is uniform in z-direction. 
The visualization of the computational mesh on the x-y plane is shown in Figure \ref{cylinder-mesh}.
The mesh is refined behind the cylinder to accurately capture the wake region, while it is coarsened in the far-field region to reduce computational resources.

In Figure \ref{cylinder-vortex}, the vortex structures around the cylinder are visualized using the iso-surface of the Q-criterion ($Q = 0.01$), colored by the Mach number, after the turbulent flow has fully developed.
The small-scale structures in the wake region are well resolved using the current fifth-order CGKS.
Figure \ref{cylinder-line} shows the averaged streamwise velocity profiles along the transverse direction at various locations in the cylinder wake.
The operator $⟨(\cdot)⟩$ represents temporal averaging and spatial averaging in the z-direction. Temporal averaging is performed at each time step after the flow achieves a statistically stationary state, free from any initial transients.
In Figure \ref{cylinder-line}, the results computed by the linear CGKS are compared with experimental data obtained using the smoke image velocimetry (SIV) technique \cite{Cylinder-ref-2}, as well as with the results from the LES method combined with a sixth-order compact scheme on a refined $256 \times 128 \times 128$ mesh \cite{Cylinder-ref-1}.
It can be observed that the results from the linear CGKS-5 closely match the experimental data, particularly in the near wake region ($x/D < 3$).
Compared to the LES results \cite{Cylinder-ref-1}, the current scheme shows a clear advantage in resolution on a mesh of the same size.

\begin{figure}[!h]
\centering
\includegraphics[width=0.72\textwidth]{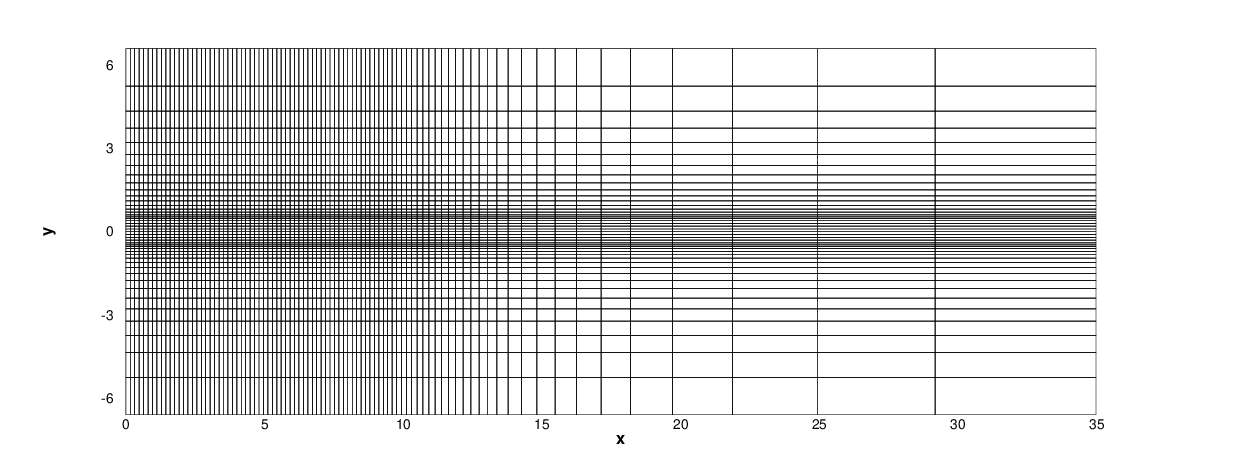}
\includegraphics[width=0.27\textwidth]{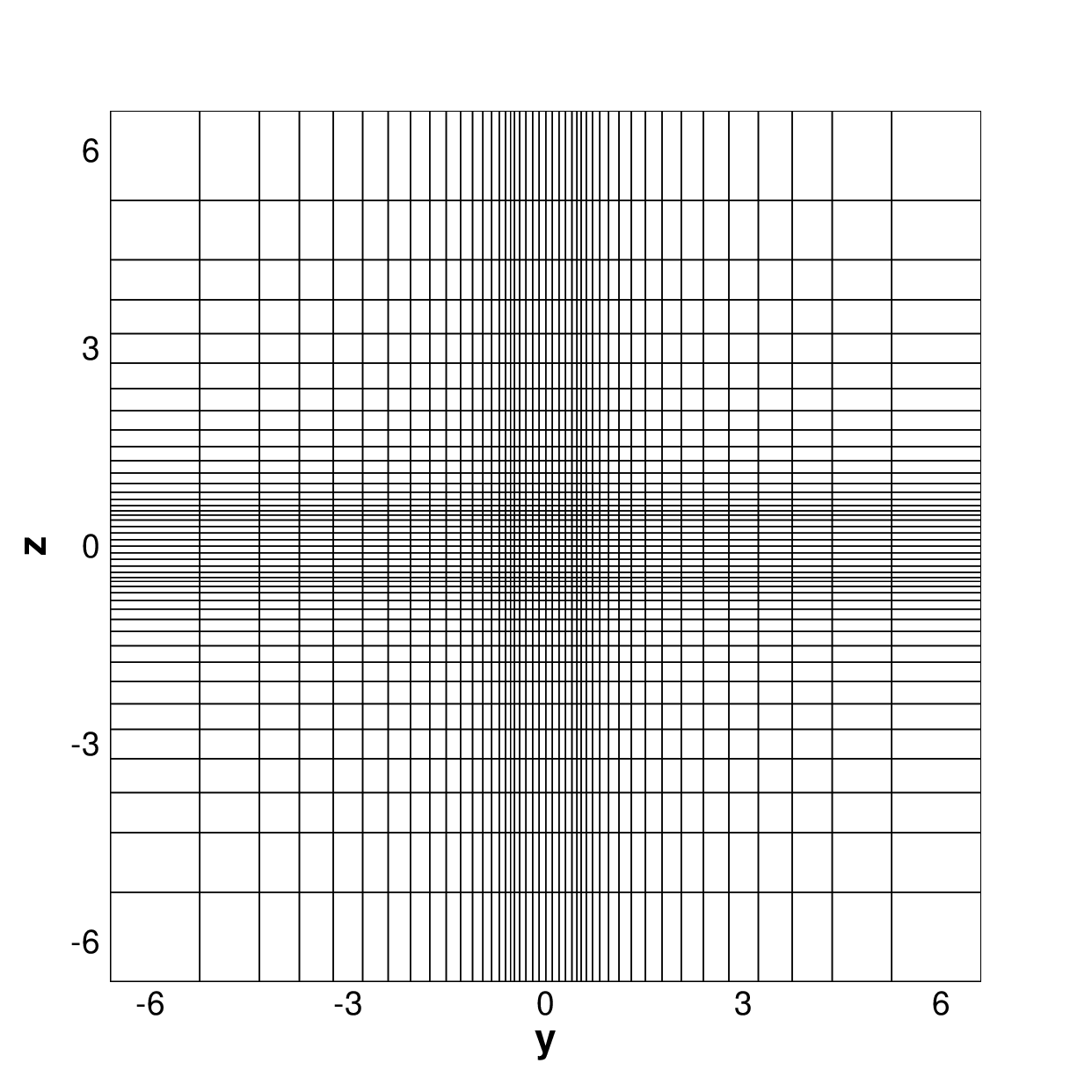}
\caption{\label{jet-mesh} Transonic and Supersonic Jet Simulation: the visualization of the computational mesh on the $y=0$ plane (left) and the $x=0$ plane (right). Only every fourth line in the three coordinate directions is shown.}
\end{figure}

\begin{figure}[!h]
\centering
\includegraphics[width=0.7\textwidth]{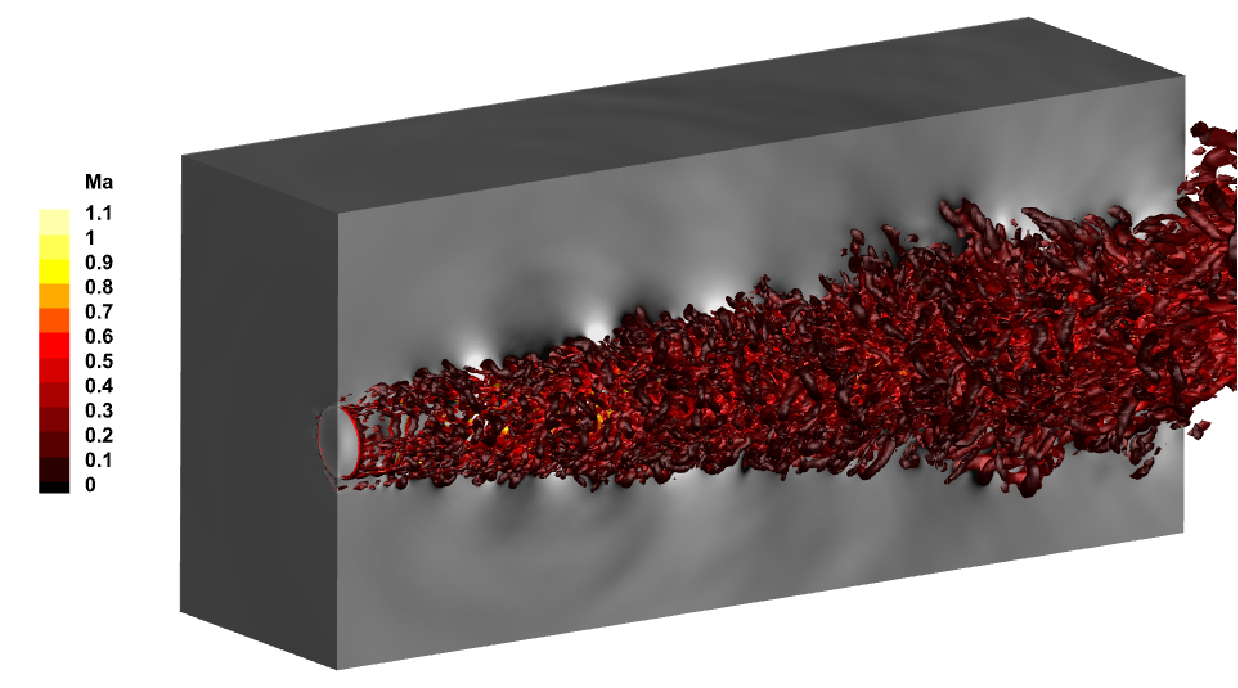}
\caption{\label{jet-sub-Q} Transonic Jet Simulation: the iso-surface of the Q-criterion ($Q=0.2$) shown in an isometric view, along with the pressure fluctuations around $P_0 \pm 10^{-2}P_0$ visualized in grayscale.}
\includegraphics[width=0.7\textwidth]{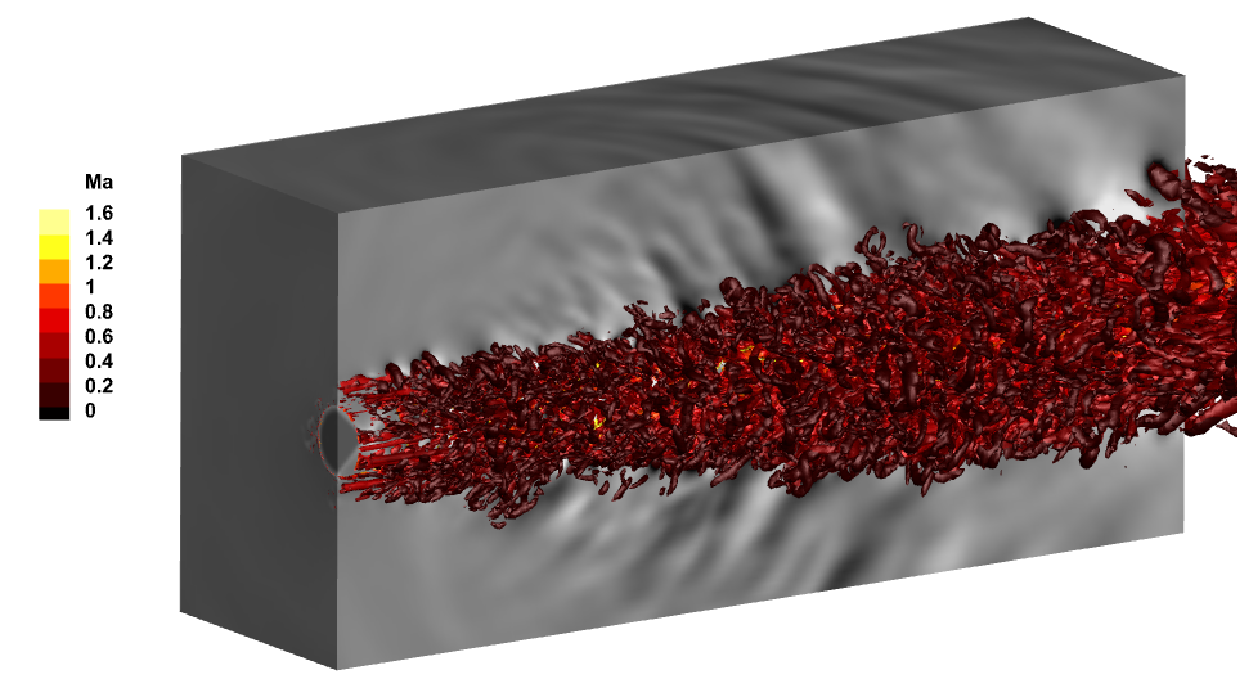}
\caption{\label{jet-sup-Q} Supersonic Jet Simulation: the iso-surface of the Q-criterion ($Q=0.4$) shown in an isometric view, along with the pressure fluctuations around $P_0 \pm 4\times 10^{-2}P_0$ visualized in grayscale.}
\end{figure}

\begin{figure}[!h]
\centering
\includegraphics[width=0.485\textwidth]{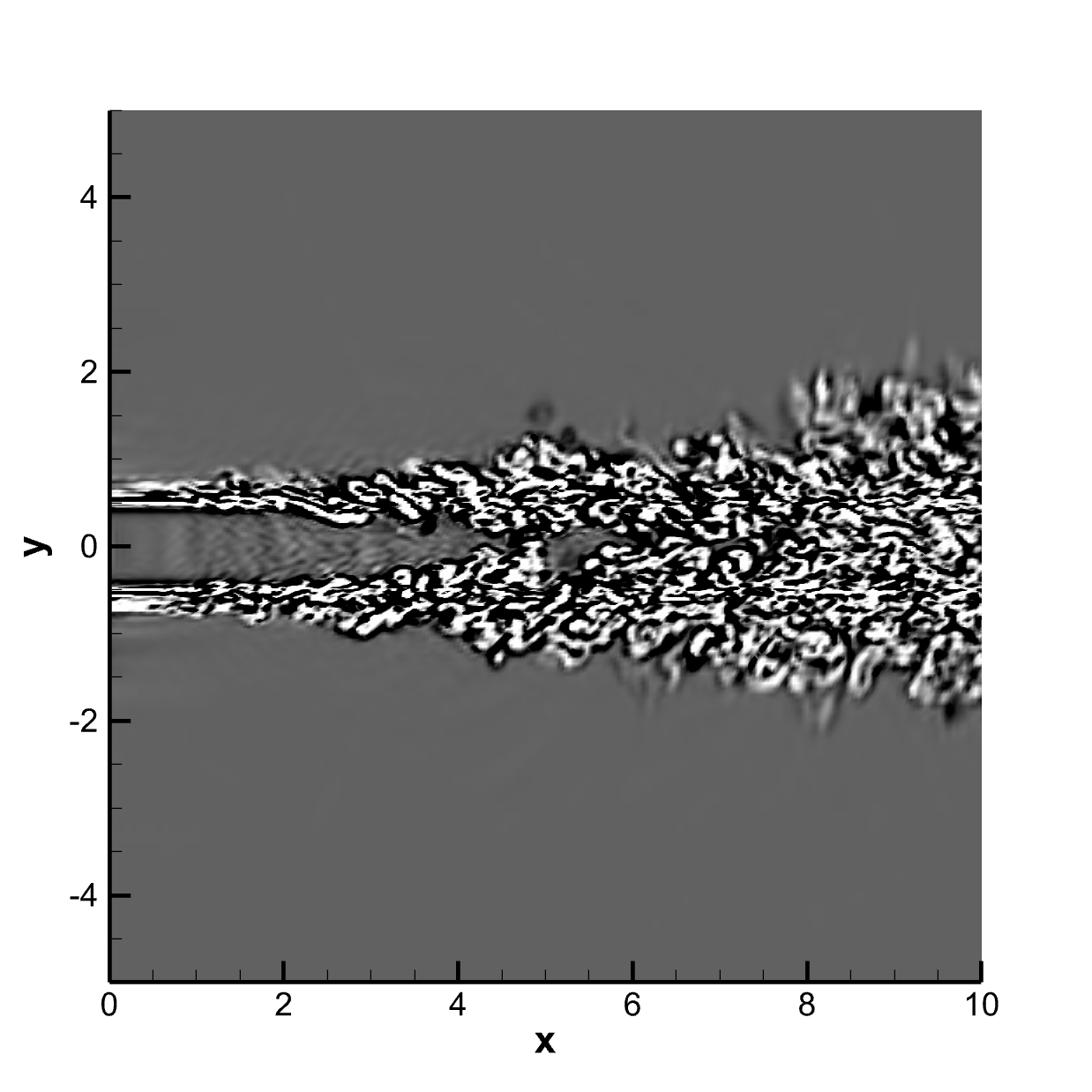}
\includegraphics[width=0.485\textwidth]{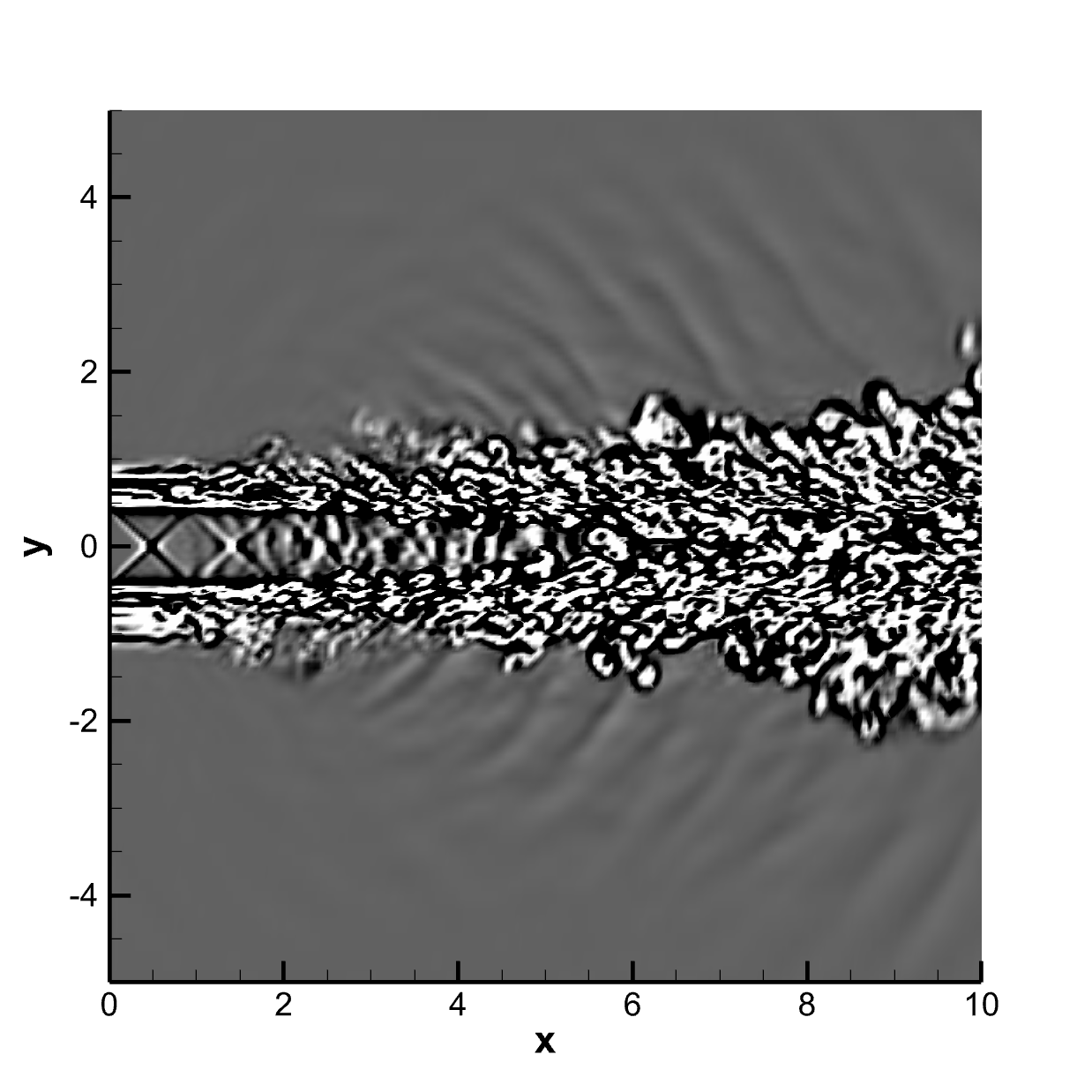}
\caption{\label{jet-pre} Transonic and Supersonic Jet Simulation: the instantaneous shadowgraph for $Ma_{\infty}=0.9$ (left) and $Ma_{\infty}=1.5$ (right).}
\end{figure}

\begin{figure}[!h]
\centering
\includegraphics[width=0.485\textwidth]{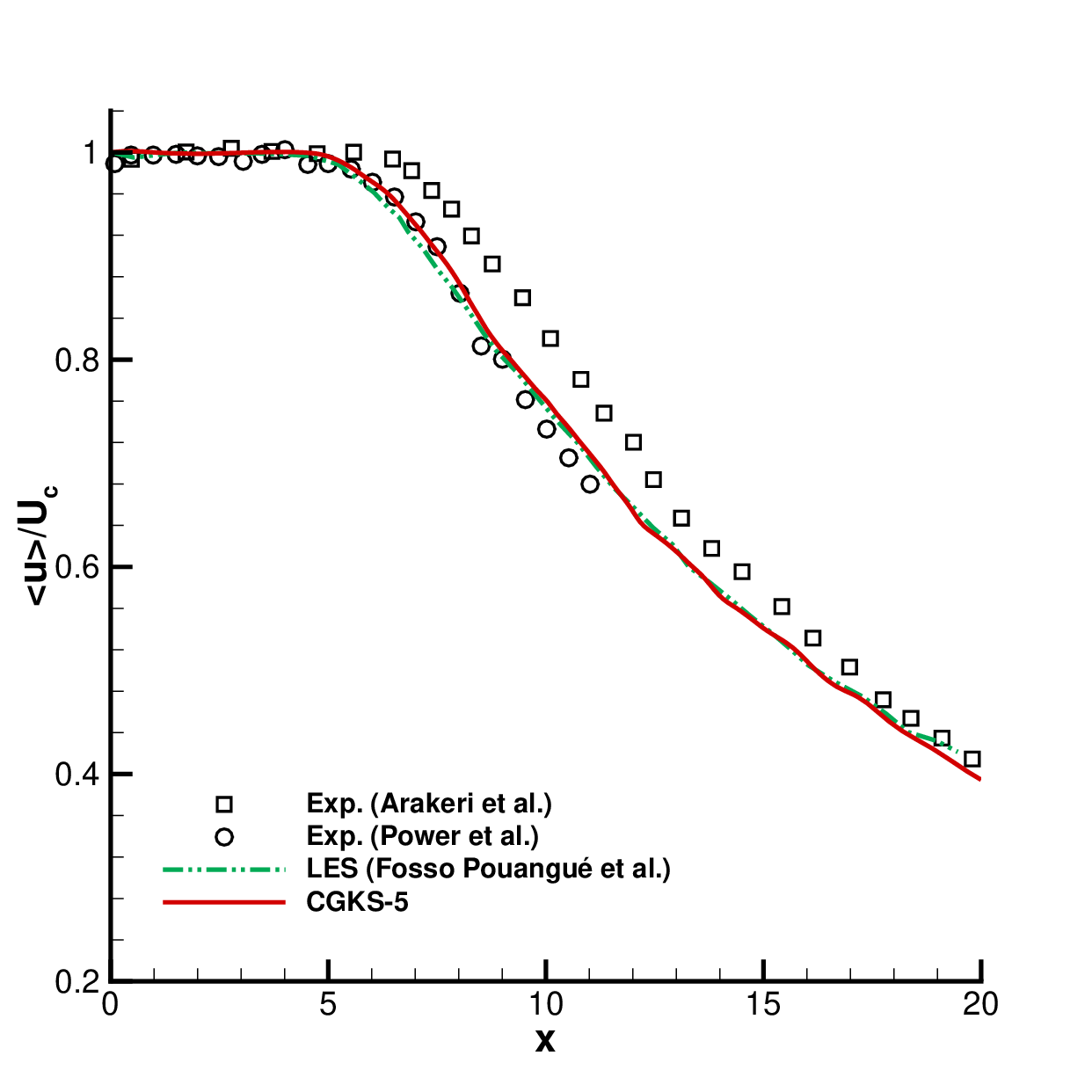}
\includegraphics[width=0.485\textwidth]{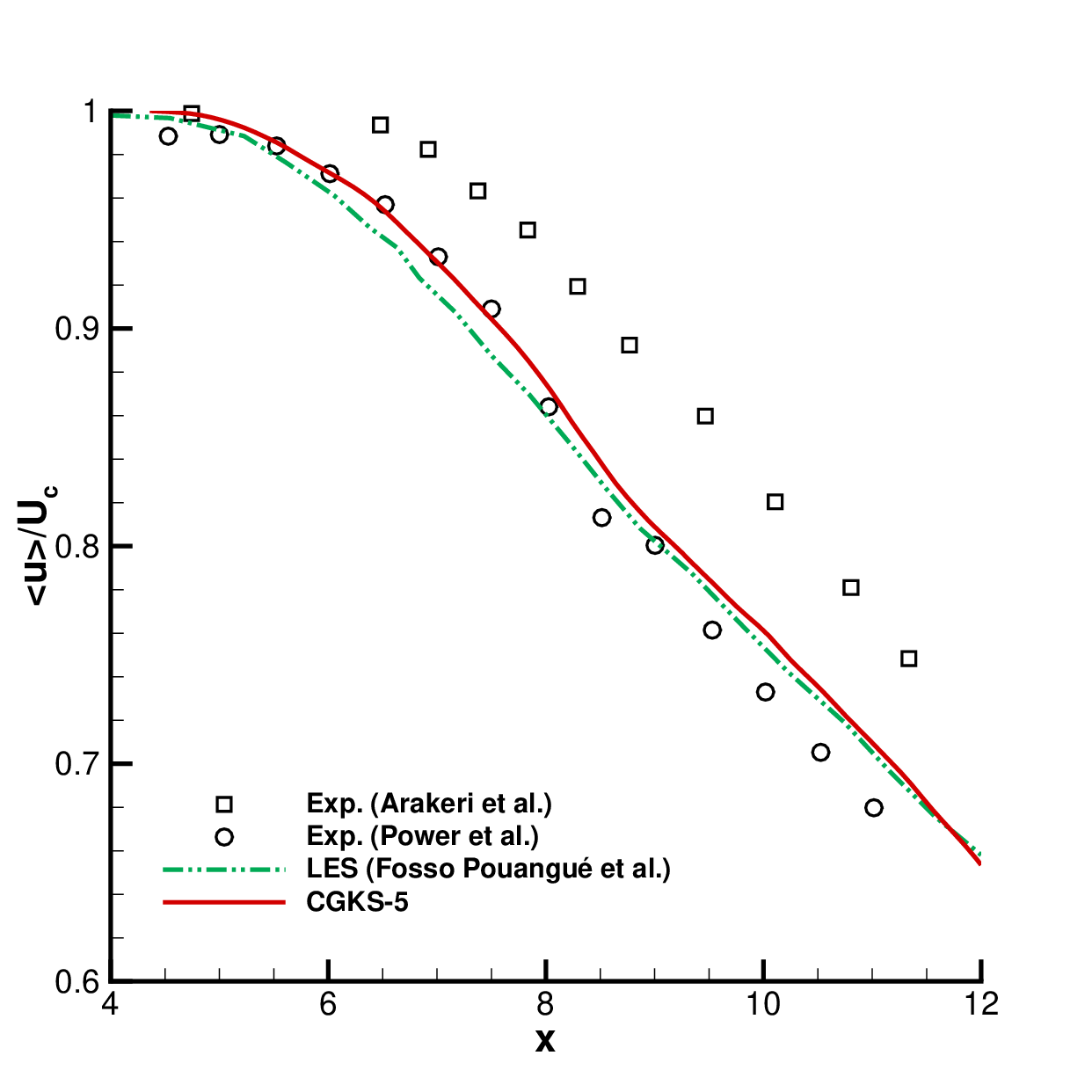}
\caption{\label{jet-sub-ave1} Transonic Jet Simulation: the mean axial velocity (left) and local magnification (right) along the centerline.}
\includegraphics[width=0.485\textwidth]{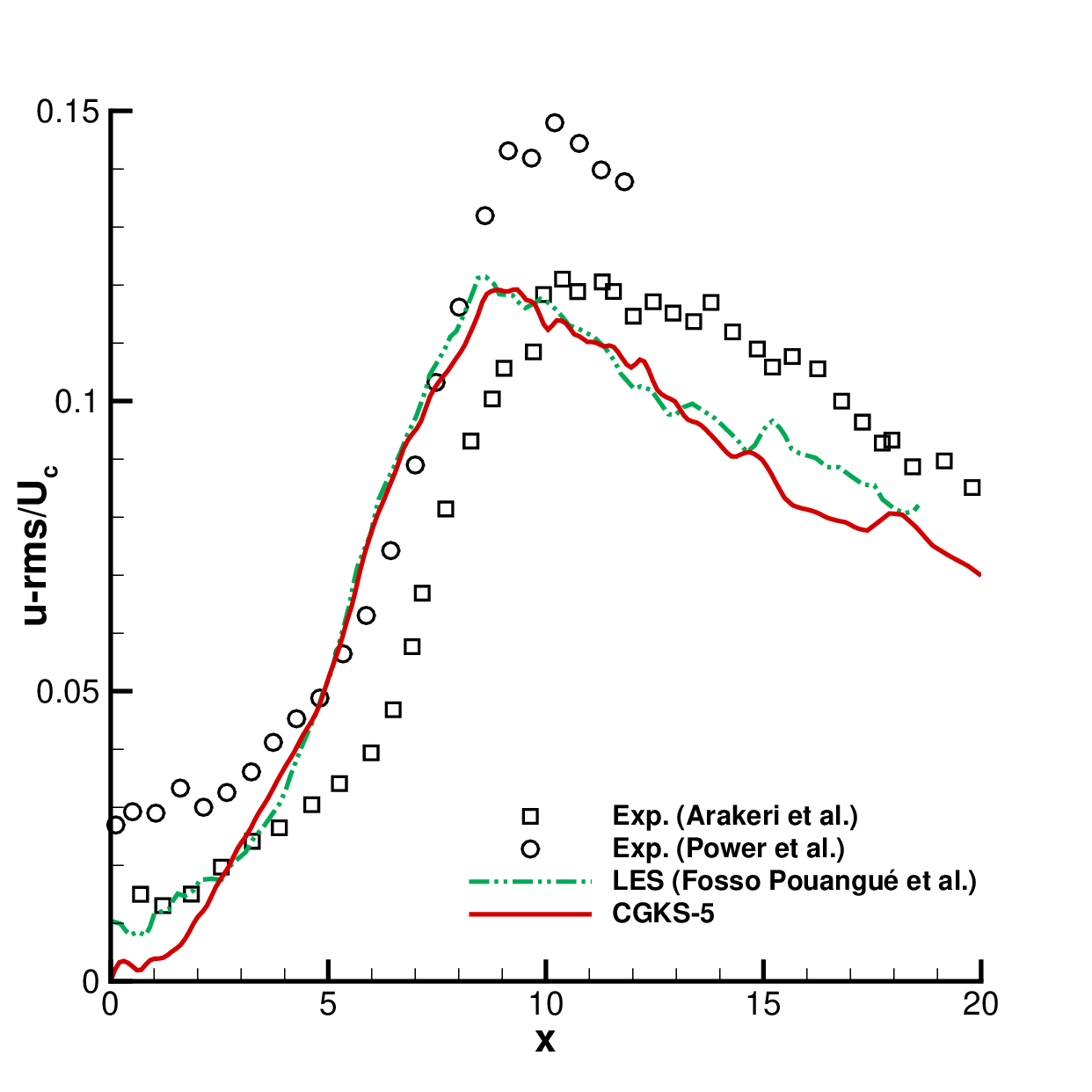}
\includegraphics[width=0.485\textwidth]{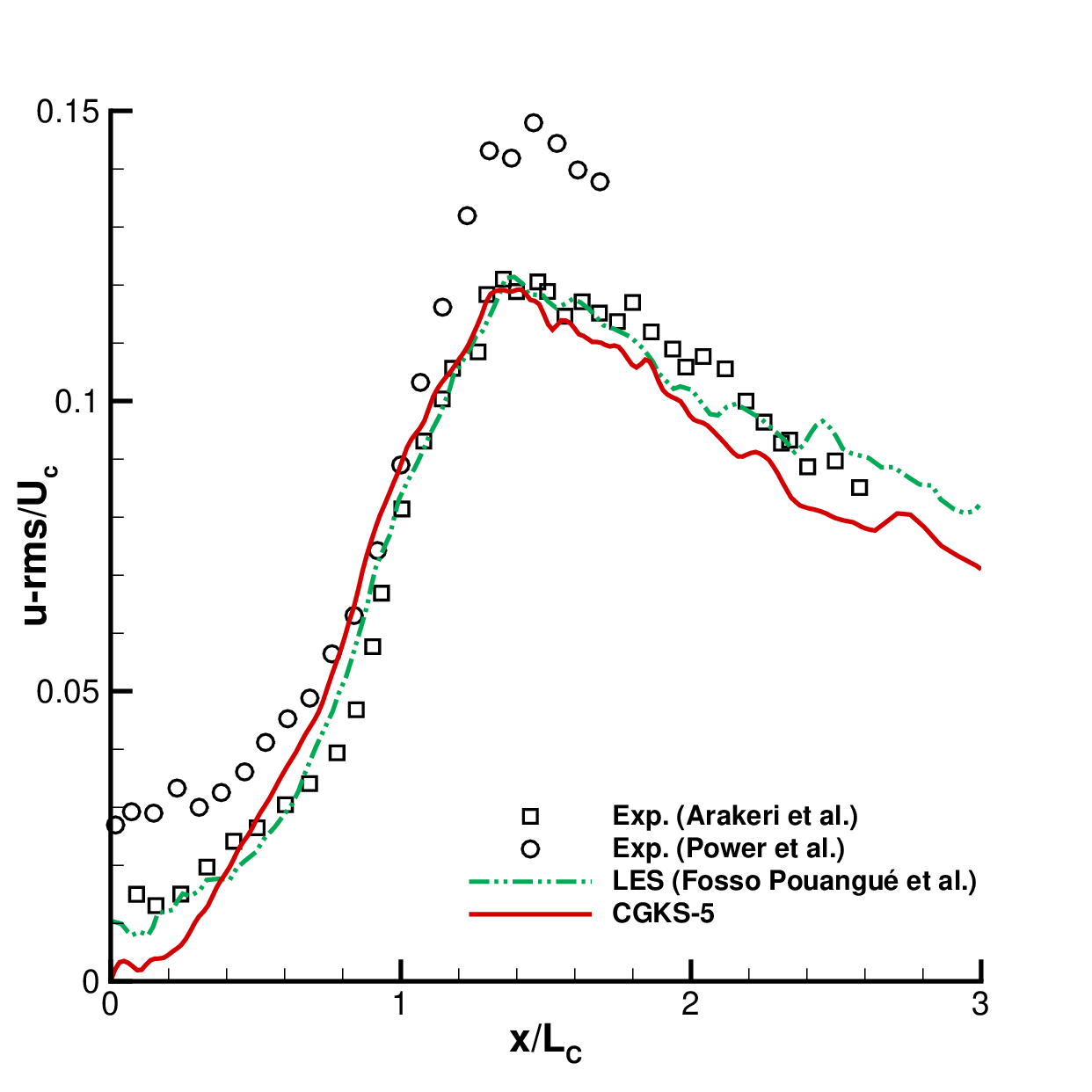}
\caption{\label{jet-sub-ave2} Transonic Jet Simulation: the root mean square axial velocity (left) and $L_C$-normalized distribution (right) along the centerline.}
\end{figure}

\begin{figure}[!h]
\centering
\includegraphics[width=0.485\textwidth]{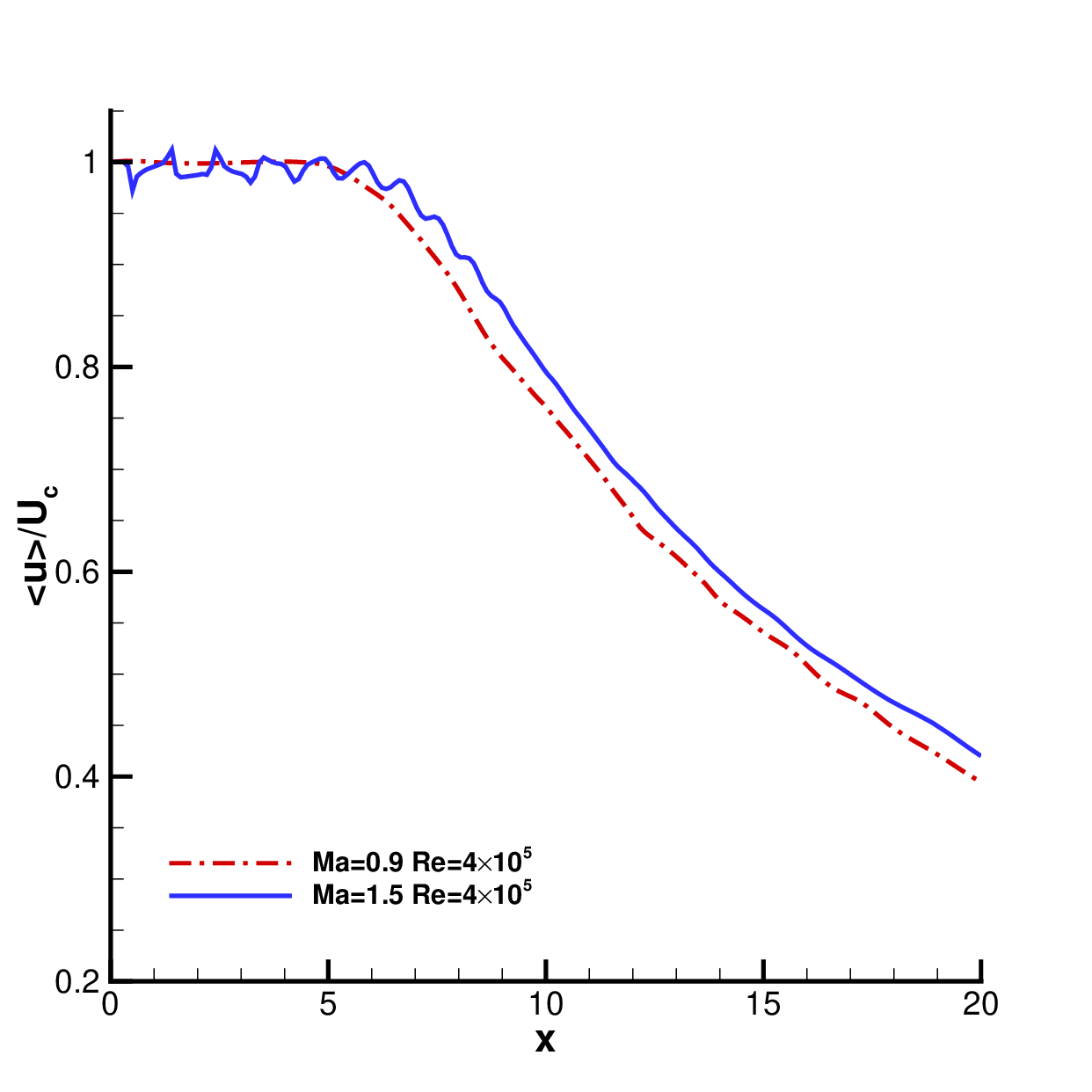}
\includegraphics[width=0.485\textwidth]{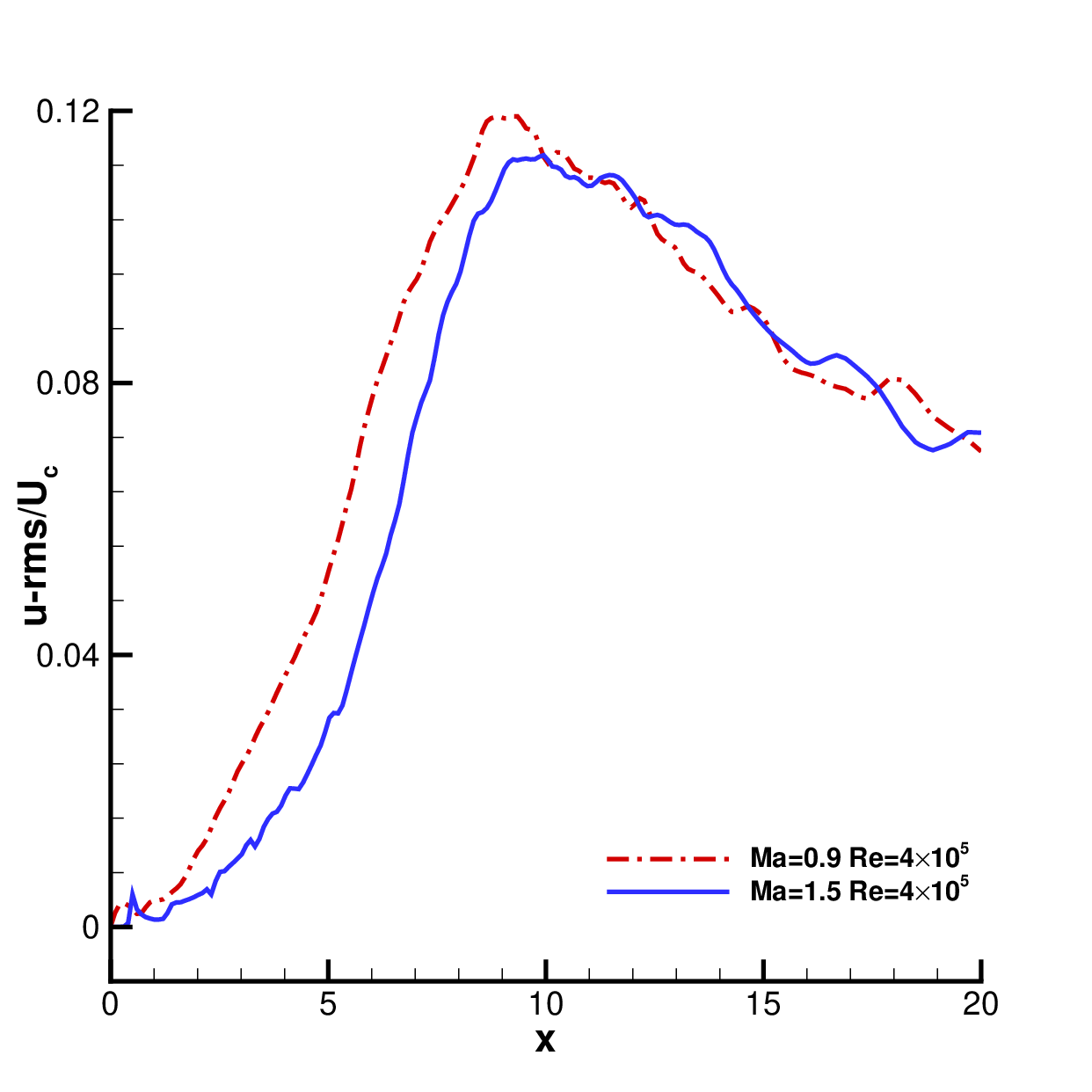}
\caption{\label{jet-sup-ave} Supersonic Jet Simulation: the mean axial velocity (left) and root mean square axial velocity (right) along the centerline.}
\end{figure}

\subsection{Transonic and Supersonic Jet Simulation}

Jet simulation plays an important role in understanding and optimizing jet flows for a variety of engineering applications.
In the simulation, jets typically exhibit high-velocity gradients at the boundary between the jet and the surrounding fluid. 
As jets develop downstream, they generally transition from laminar to turbulent flow.
In this subsection, Transonic and supersonic jet simulations are selected to validate the robustness and resolution of the current fifth-order compact GKS.
The isothermal jets with Mach numbers $Ma_{\infty} = 0.9$ and $Ma_{\infty} = 1.5$ are considered, with the Reynolds number set as $Re=4\times 10^5$ based on a jet diameter of $D =2r_0= 1$.
In this subsection, None of the simulations include any nozzle geometry.

In this case, a Cartesian mesh consisting of $352\times 176\times 176$ cells is stretched, as illustrated in Figure \ref{jet-mesh}, to meet computational requirements.
In the flow direction, the mesh size is constant up to $x=9.6D$, with $\Delta x=0.04D$, to accurately resolve the physical domain near the nozzle exit.
A transition region $9.6DD<x<15D$ is defined, where the mesh size increases exponentially from $0.04D$ to $0.16D$.
Then the mesh is stretched to create a sponge zone for damping the vortices and perturbation, where $\Delta x_{max}=1.6D$.
In the radial direction, the mesh size at the jet center is $0.025D$. To resolve the small-scale structures at the boundary between the jet and the surrounding fluid near the jet, the mesh size decreases exponentially to a minimum of approximately $0.015D$ at $r=r_0$. 
Outside the jet, the mesh size increases back to $0.025D$, then gradually expands to $0.2D$ within $r<4.5D$. Beyond this region, a sponge region is formed with the maximum mesh size reaching $0.4D$.

The inflow axial velocity of the jet is initialized through a hyperbolic-tangent profile similarly to previous studies \cite{Jet-ref-1, Jet-ref-2}
\begin{equation*}
u(r)=\displaystyle\frac{U_c}{2}+\frac{U_c}{2}\tanh (\frac{r_0-r}{2\delta_\theta(0)}),
\end{equation*}
where $U_c=Ma_{\infty}\cdot c_\infty$ is the inflow centerline velocity, $c_\infty=1$, and the initial momentum thickness $\delta_\theta(0)$ is set as $0.03D$.
The inlet temperature profile is calculated using the Crocco-Busemann relation, assuming constant pressure
\begin{equation*}
\displaystyle\frac{T}{T_c}=\frac{T_\infty}{T_c}+\big[1-\frac{T_\infty}{T_c}+\frac{\gamma-1}{2}{Ma_{\infty}}^2\big(1-\sigma \frac{u(r)}{U_c}\big)\big]\frac{u(r)}{U_c},
\end{equation*}
where $T_c$ and $T_\infty$ are the jet centerline and ambient temperatures, respectively, in the present isothermal case, $T_c=T_\infty$; and $\sigma={Pr}^{1/3}$ is the recovery factor, the Prandtl number $Pr$ is set to 0.7.
Outside the inlet, the initial condition is given as follows,
\begin{align*}
(\rho,U,V,W,p)_c = (\gamma, 0, 0, 0, 1).
\end{align*}
No disturbance is added to the field.
The non-reflecting boundary condition and sponge zone are implemented in the outflow directions to filter out possible reflected waves.
For the transonic jet, large-scale vortex ring structures dominate during the development stage. Due to the high Reynolds numbers, disturbances within the shear layer are rapidly amplified. Interactions between these disturbances and shock waves result in a rapid transition to turbulence structure. 
The supersonic jet exhibits stronger shock-shear layer interactions, leading to a faster transition and a more complex turbulence structure.

Figures \ref{jet-sub-Q} and \ref{jet-sup-Q} show instantaneous iso-surfaces of the Q-criterion, colored by the Mach number, for $Ma_{\infty}=0.9$ and $Ma_{\infty}=1.5$, respectively. The complex multi-scale turbulent structures of the jets are clearly visible. Additionally, Figures \ref{jet-sub-Q} and \ref{jet-sup-Q} include instantaneous pressure grayscale images depicting pressure fluctuations relative to $P_0$, while Figure \ref{jet-pre} shows the corresponding instantaneous shadowgraph contours.
These visualizations highlight the propagation of sound waves within the computational domain, with the supersonic jet ($Ma_{\infty}=1.5$) demonstrating more pronounced and distinct sound waves compared to the subsonic jet.

Figure \ref{jet-sub-ave1} and Figure \ref{jet-sub-ave2} illustrate the mean axial velocity and the root mean square (rms) axial velocity along the centerline ($\boldsymbol{r}=0$) with $Ma_{\infty}=0.9$, respectively.
In Figure \ref{jet-sub-ave1} and Figure \ref{jet-sub-ave2}, the numerical results obtained by CGKS-5 are compared to the reference results, including experimental results from Arakeri et al.\cite{Jet-ref-4} and Power et al.\cite{Jet-ref-3}, as well as LES result using a sixth-order compact scheme on a multiblock structured mesh with 24.5 million cells (Fosso Pouangu{\'e} et al.)\cite{Jet-ref-2}.
In Figure \ref{jet-sub-ave2}, the axial position is scaled by the corresponding potential core length, $L_C$, to facilitate a comparison of the turbulent intensity development relative to the decay of the jet centerline velocity.
The potential core length is defined as the distance from the jet inlet to the point where the centerline velocity reduces to 95$\%$ of the inflow centerline velocity. 
Table \ref{jet-Lc} summarizes the $L_C$ computed by CGKS-5 and compares it with reference data. It can be observed that the $L_C$ calculated using CGKS-5 show good agreement with the reference data.
% In the computation, the potential core length $L_C$ is about $6.6D$ for the transonic jet and $7.1D$ for the supersonic jet.
In Figures \ref{jet-sub-ave1} and \ref{jet-sub-ave2}, the mean axial velocity results computed using CGKS-5 exhibit excellent agreement with the experimental data as well as with the LES result obtained using finer mesh and higher-order discretization. Similarly, the rms axial velocity results calculated by CGKS-5 align closely with the LES numerical results \cite{Jet-ref-2} and the experimental data from Arakeri et al.\cite{Jet-ref-4}. These findings demonstrate the high accuracy and resolution of CGKS-5 in the scale-resolving simulation of complex turbulent flows.
% Furthermore, all results indicate that the maximum turbulence intensity occurs around $x=1.4L_C$.

\begin{table}[!h]
\begin{center}
\def\temptablewidth{0.5\textwidth}{\rule{\temptablewidth}{1.0pt}}
\begin{tabular*}{\temptablewidth}{@{\extracolsep{\fill}}c|c} 
Method & $L_C$  \\
\hline
Exp. Arakeri et al. \cite{Jet-ref-4} & 7.7D \\
Exp. Power et al. \cite{Jet-ref-3} & 7.0D \\
LES Fosso Pouangu{\'e} et al. \cite{Jet-ref-2} & 6.2D \\
CGKS-5 & 6.6D 
\end{tabular*}
{\rule{\temptablewidth}{1.0pt}}
\end{center}
\caption{\label{jet-Lc} Transonic Jet Simulation: the comparison of the $L_C$ calculated by CGKS-5 with reference data.}
\end{table}

Furthermore, Figure \ref{jet-sup-ave} presents the mean axial velocity and rms axial velocity for different Mach numbers at the same Reynolds number. It can be observed that as the Mach number increases, the potential core length extends, while the turbulence intensity is suppressed due to the increased viscosity.
Notably, in Figure \ref{jet-sup-ave}, the fluctuations in the mean axial velocity result obtained by CGKS-5 within the region $x \in [0,6]$ correspond to the locations of shock waves inside the jet observed in Figure \ref{jet-pre}.

\section{Conclusion}

In this paper, we present a novel high-resolution CGKS for three-dimensional structured meshes.
A new multidimensional fifth-order compact reconstruction method is introduced for Cartesian meshes, utilizing line-averaged derivatives within each control volume to construct higher-order stencils and achieve enhanced resolution.
To extend the new CGKS to non-orthogonal meshes, a straightforward coordinate transformation approach is employed. 
Reconstruction is then performed on the transformed computational mesh, enabling higher computational efficiency and significantly reduced memory requirements. 
Additionally, an efficient nonlinear combination mechanism is employed to ensure both high resolution and robustness.
Furthermore, the proposed method is implemented on multiple GPUs using CUDA, with numerical tests conducted on 8 NVIDIA GeForce RTX 4090 GPUs.
The results demonstrate excellent resolution, stability, and robustness of the new fifth-order CGKS in simulating turbulent flows, making it a promising tool for high-fidelity, large-scale simulations of complex three-dimensional turbulence.

\section*{Acknowledgements}

The current research is supported by National Science Foundation of China (92371107, 12172316), National Key R$\&$D Program of China (Grant Nos.
2022YFA1004500), and Hong Kong research grant council (16301222, 16208324).

\section*{Declaration of competing interest}

The authors declare that they have no known competing financial interests or personal relationships 
that could have appeared to influence the work reported in this paper.

\section*{Data availability}

The data that support the findings of this study are available from
the corresponding author upon reasonable request.

\section*{Appendix}

A trilinear interpolation is introduced for the transformed neighboring cell $\Omega^{'}\in N_f(0)^{'}\cup N_e(0)^{'}$
\begin{align*}
\hat{\boldsymbol{X}}(\xi,\eta,\zeta)=\sum_{m=1}^8\boldsymbol{x}_m\psi_m(\xi,\eta,\zeta),
\end{align*}
where $(\xi,\eta,\zeta)\in[-1/2,1/2]^3$, $\boldsymbol{x}_m$ is the
vertex of cell $\Omega^{'}$ and the basis function $\psi_m$ is given
as follows
\begin{align*}
\psi_1&=\frac{1}{8}(1-2\xi)(1-2\eta)(1-2\zeta),\quad
\psi_2=\frac{1}{8}(1-2\xi)(1-2\eta)(1+2\zeta),\\
\psi_3&=\frac{1}{8}(1-2\xi)(1+2\eta)(1-2\zeta),\quad
\psi_4=\frac{1}{8}(1-2\xi)(1+2\eta)(1+2\zeta),\\
\psi_5&=\frac{1}{8}(1+2\xi)(1-2\eta)(1-2\zeta),\quad
\psi_6=\frac{1}{8}(1+2\xi)(1-2\eta)(1+2\zeta),\\
\psi_7&=\frac{1}{8}(1+2\xi)(1+2\eta)(1-2\zeta),\quad
\psi_8=\frac{1}{8}(1+2\xi)(1+2\eta)(1+2\zeta).
\end{align*}
In the computation, the Gaussian quadrature is used for the
following triple integral over the cell $\Omega^{'}$
\begin{align*}
|\Omega^{'}|=\int_{\Omega^{'}}1\,\text{d}\hat{x}\text{d}\hat{y}\text{d}\hat{z}
=\sum_{l,m,n=1}^3\omega_{lmn}\Big|\frac{\partial(\hat{x},\hat{y},\hat{z})}{\partial(\xi,\eta,\zeta)}\Big|_{(\xi_l,\eta_m,\zeta_n)}\Delta\xi\Delta\eta\Delta\zeta,
\end{align*}
where $\omega_{lmn}$ is the Gaussian quadrature weight and
$(\xi_l,\eta_m,\zeta_n)$ is the quadrature point.

\end{document}